%&Plain
%% macroarticolops.tex 
%% per la data di sistema
\def\oggi{\number\day\space\ifcase\month\or
       gennaio\or febbraio\or marzo\or aprile\or maggio\or giugno\or
       luglio\or agosto\or settembre\or ottobre\or novembre\or dicembre\fi
        \ \number\year}
\def\today{\ifcase\month\or
               January\or February\or March\or April\or May\or June\or
               July\or August\or September\or October\or November\or December\fi
               \space\number\day, \number\year}
\def\aujourdhui{\number\day\space\ifcase\month\or
               Janvier\or F\'evrier\or Mars\or Avril\or May\or Juin\or
               Juillet\or Ao\^ut\or Septembre\or Octobre\or Novembre\or 
               D\'ecembre\fi
               \ \number\year}

%%%%%%%%%%%%%%%%%%%%%%%%%%%%%%%%%%%%%%%%%%%%%%%%%%%%%%%%%%%%%%%%%%%%%%%%%%
%%%%%%%%%%%%%%% contatori fino a quattro entrate
%%%%%%%%%%%%%%% uso: \N o \NN o \NNN o \NNNN
%%%%%%%%%%%%%%% effetto: stampa numerazione nel font corrente
%%%%%%%%%%%%%%%%%%%%%%%%%%%%%%%%%%%%%%%%%%%%%%%%%%%%%%%%%%%%%%%%%%%%%%%%%%

\newcount\s
\newcount\n
\def\clearn{\n =0}
\def\N{\global\advance\n by 1 \global\s=1
                       \global\clearnn
                       \global\clearnnn
                       \global\clearnnnn
                       {\the\n}%
                     }
\newcount\nn
\def\clearnn{\nn =0}
\def\NN{\global\advance\nn by 1 \global\s=2
                        \global\clearnnn
                        \global\clearnnnn
                        {\the\n}.{\the\nn}%
                   }
\newcount\nnn
\def\clearnnn{\nnn =0}
\def\NNN{\global\advance\nnn by 1 \global\s=3
                        \global\clearnnnn
                        {\the\n}.{\the\nn}.{\the\nnn}%
                     }

\newcount\nnnn
\def\clearnnnn{\nnnn =0}
\def\NNNN{\global\advance\nnnn by 1 \global\s=4
                        {\the\n}.{\the\nn}.{\the\nnn}.{\the\nnnn}%
                     }

%%%%%%%%%%%%%%%%%%%%%%%%%%%%%%%%%%%%%%
%%%%%%%%%%%%%%%%%%%%%%%%%%%%%%%%%%%%%%
%%%%%%%%%%%%%%%%%%%%%%%%%%%%%%%%%%%%%%
%%%%%%%%%%%%%%%%%%%%%%%%%%%%%%%%%%%%%%%%%%%%%%%%%%%%%%%%%%%%%%%%%%%%%
\newif\ifindi
\newbox\boxindice
\inditrue
%\indifalse
\newwrite\ind
\def\apriindice{\ifindi%
                            \setbox\boxindice=\vbox{\input \nomefile.ind }%
                            %\long\edef\copiaindice{\boxindice }%
                            \immediate\openout\ind=\nomefile.ind%
                            \immediate\write\ind{{\sectionfont
Contents}\vskip10pt}%
                            \else\fi}
\def\indice#1{\ifindi\immediate\write\ind{{#1}}\else\fi}

%\def\stampaindice{\ifindi\copiaindice\else\fi}
%%%%%%%%%%%%%%%%%%%%%%%%%%%%%%%%%%%%%%%%%%%%%%%%%%%%%%%%%%%%%%%%%%%%%
%%%%%%%%%%%%%%%%%%%%%%%%%%%%%%%%%%%%%%
%%%%%%%%%%%%%%%%%%%%%%%%%%%%%%%%%%%%%%
%%%%%%%%%%%%%%%%%%%%%%%%%%%%%%%%%%%%%%

%%%%%%%%%%%%%%%%%%%%%%%%%%%%%%%%%%%%%%%%%%%%%%%%%%%%%%%%%%%%%%%%%%
\newif\ifrifer
\rifertrue
%\riffalse
\newwrite\rifer
\def\apririferimenti{\ifrifer%
                            \immediate\openout\rifer=\nomefile.rif%
                            \immediate\write\rifer{riferimenti}%
                            \else\fi}
\def\riferimenti#1{\ifrifer\immediate\write\rifer{#1}\else\fi}
\def\chiudiriferimenti{\ifrifer\immediate\closeout\rifer\else\fi}
%%%%%%%%%%%%%%%%%%%%%%%%%%%%%%%%%%%%%%%%%%%%%%%%%%%%%%%%%%%%%%%%%%

%%%%%%%%%%%%%%%%%%%%%%%%%%%%%%%%%%%%%%%%%%%%%%%%%%%%%%%%%%%%%%%%%%%%%
%%%%%%%%%%%%%% \rif{NOME} crea le label per il riferimento
%%%%%%%%%%%%%%  (lo stato della numerazione in quel momento)
%%%%%%%%%%%%%% \cite{NOME} scrive la citazione corrispondente a NOME
%%%%%%%%%%%%%%  (nel font corrente)
%%%%%%%%%%%%%%%%%%%%%%%%%%%%%%%%%%%%%%%%%%%%%%%%%%%%%%%%%%%%%%%%%%%%%

\def\rif#1{\expandafter\xdef\csname s#1\endcsname{\number\s}%
            \ifnum\s>0 \expandafter\xdef\csname n#1\endcsname{\number\n}\fi%
            \ifnum\s>1 \expandafter\xdef\csname nn#1\endcsname{\number\nn}\fi%
            \ifnum\s>2 \expandafter\xdef\csname nnn#1\endcsname{\number\nnn}\fi%
            \ifnum\s>3 \expandafter\xdef\csname nnnn#1\endcsname{\number\nnnn}\fi%
            %\riferimenti{\the\n.\the\nn.\the\nnn.\the\nnnn.#1}%
            \riferimenti{\cite{#1}........#1}%
                   }

\def\cite#1{%
            \expandafter\ifcase\csname s#1\endcsname%
            \or \csname n#1\endcsname%
            \or \csname n#1\endcsname.\csname nn#1\endcsname%
            \or \csname n#1\endcsname.\csname nn#1\endcsname.\csname nnn#1\endcsname%
            \or \csname n#1\endcsname.\csname nn#1\endcsname.\csname nnn#1\endcsname.%
                \csname nnnn#1\endcsname%
                   \fi}

%%%%%%%%%%%%%%%%%%%%%%%%%%%%%%%%%%%%%%%%%%%%%%%%%%%%%%%%%%%%%%%%%%%%%%%%%
\font \proc=cmbx10
\def\proclaim#1. #2\par
                  {\medbreak  {\proc #1. \enspace }{\sl #2\par }
                    \ifdim \lastskip <\medskipamount \removelastskip
                            \penalty 55\medskip \fi}
%%%%%%%%%%%%%%%%%%%%%%%%%%%%%%%%%%%%%%%%%%%%%%%%%%%%%%%%%%%%%%%%%%%%%%%%%

%%%%%%%%%%%%%%%%%%%%%%%%%%%%%%%%%%%%%%%%%%%%%%%%%%%%%%%%%%%%%%%%%%%%%%%%%
\def\biblio#1#2{ \item{\hbox to 1.2cm{[#1]\hss}}
                           {#2} }
%%%%%%%%%%%%%%%%%%%%%%%%%%%%%%%%%%%%%%%%%%%%%%%%%%%%%%%%%%%%%%%%%%%%%%%%%

%%%%%%%%%%%%%%%%%%%%%%%%%%%%%%%%%%%%%%%%%%%%%%%%%%%%%%%%%%%%%%%%%%%%%%%%%
%% ambienti predefiniti #1=numerazione (N, NN, NNN o NNNN)
%%                      #2=testo

\def\numfont{\bf}
%\font\numfont=?

\def\prop#1#2{\proclaim{{\numfont\csname#1\endcsname.} Proposition}. {#2}\par}
\def\theo#1#2{\proclaim{{\numfont\csname#1\endcsname.} Theorem}. {#2}\par}
\def\lemm#1#2{\proclaim{{\numfont\csname#1\endcsname.} Lemma}. {#2}\par}
\def\defi#1#2{\proclaim{{\numfont\csname#1\endcsname.} Definition}. {#2}\par}
\def\coro#1#2{\proclaim{{\numfont\csname#1\endcsname.} Corollary}. {#2}\par}

\long\def\proo#1#2{{{\numfont\csname#1\endcsname.}\proc \ Proof.}
                                    #2
                                    \hfill$\square$\vskip 10pt
                                }
\long\def\exam#1#2{{{\numfont\csname#1\endcsname.}\proc \ Example.}
                                    #2
                                    \vskip 10pt
                                }

\long\def\rema#1#2{{{\numfont\csname#1\endcsname.}\proc \ Remark.}
                                    #2
                                    \vskip 10pt
                                }

\long\def\nota#1#2{{{\numfont\csname#1\endcsname.}\proc \ Notations.}
                                    #2
                                    \vskip 10pt
                                }

\long\def\proof#1{{\proc \ Proof.}
                                    #1
                                    \hfill$\square$\vskip 10pt
                                }

%%%%%%%%%%%%%%%%%%%%%%%%%%%%%%%%%%%%%%%%%%%%%%%%%%%%%%%%%%%%%%%%%%%%%%%%%
%% ambienti predefiniti #1=titolo

\font\titlefont=cmbx10 scaled \magstep2
\font\authorfont=cmcsc10 scaled \magstep1
\font\sectionfont=cmbx10 scaled \magstep1
\font\subsectionfont=cmbx10
\font\abstractfont=cmr8
\font\titleabstractfont=cmcsc8

\def\section#1{\vskip 10pt
                          {\sectionfont\N.\ #1}
                          \indice{\proc \the\n.\ #1\endgraf}
                          \vskip 10pt}
\def\subsection#1{\vskip 8pt
                             {\subsectionfont\NN.\ #1}
                             \indice{\quad\proc \the\nn.\ #1\endgraf}
                             \vskip 8pt}
\def\intro{\vskip 10pt
                          {\sectionfont Introduction.}
                          \vskip 10pt\indice{\proc Introduction.\endgraf}}
\def\refer{\vskip 10pt
                          {\sectionfont References.}
                          \vskip 10pt\indice{\proc References.\endgraf}}
\def\title#1{\vskip 5pt
                          \centerline{\titlefont #1}
                          \vskip 10pt}
\def\author#1{\vskip 8pt
                          %\centerline{by}\vskip 5pt
                          \centerline{\authorfont #1}
                          \vskip 10pt}

\def\abstract#1{
%\vskip 8pt
\centerline{\vtop{\hsize=12truecm\baselineskip=9pt\strut
                 {\titleabstractfont Abstract. }
                         \abstractfont #1}}
                          \vskip 10pt}

\def\resume#1{\vskip 8pt
\centerline{\vtop{\hsize=12truecm\baselineskip=9pt\strut
                 {\titleabstractfont R\'esum\'e. }
                         \abstractfont #1}}
                         % \vskip 10pt
                         }

%%%%%%%%%%%%%%%%%%%%%%%%%%%%% start %%%%%%%%%%%%%%%%%%%%%%%%%%%%%%
\clearn

%%%%%%%%%%%%%%%%%%%%%%%%%%%%%%%%%%%%%%%%% DIMENSIONS
\vsize= 215truemm
\hsize= 130truemm
\voffset= 10truemm
\hoffset= 16truemm

\abovedisplayskip=6pt plus2pt minus 4pt
\belowdisplayskip=6pt plus2pt minus 4pt
\abovedisplayshortskip=0pt plus2pt
\belowdisplayshortskip=2pt plus1pt minus 1pt

%\vsize= 225truemm
%\hsize= 155truemm
%\voffset= 0truemm
%\hoffset= 2truemm
%\topskip=15pt 

%%%%%%%%%%%%%%%%%%%%%%%%%%%%%%%%%%%%%%%%%

%%%%%%%%%%%%%%%%%%%%%%%%%%%%%%%%%%%%%%%%%%%%%%%%%%%%%%%%%%

%\headline{
%\hss\small
%\nomefile.tex --- version of \today\ ---
%Status
%\hss}

%\footline{
%\hbox to0pt{\small University of Padova, Italy\hss}
%\hss\rm\folio\hss
%\hbox to0pt{\hss\small %
%\smalltt fiorot@math.unipd.it}}

%%%%%%%%%%%%%%%%%%%%%%%%%%%%%%%%%%%%%%%%%%%%%%%%%%%%%%%%%%

%%% File: Laps.TeX
\catcode`\@=11   % Let's pretend @ is a letter
%         Vertical `laps'; cf. \llap and \rlap
\long\def\ulap#1{\vbox to \z@{\vss#1}}
\long\def\dlap#1{\vbox to \z@{#1\vss}}

%         And centered horizontal and vertical `laps'
\def\xlap#1{\hbox to \z@{\hss#1\hss}}
\long\def\ylap#1{\vbox to \z@{\vss#1\vss}}

%         And a `lap' centered on its midpoint

%         And the not-long definitions
\def\ulap#1{\vbox to \z@{\vss#1}}
\def\dlap#1{\vbox to \z@{#1\vss}}
\def\hlap#1{\hbox to \z@{\hss#1\hss}}
\def\vlap#1{\vbox to \z@{\vss\hbox{#1}\vss}}
\def\clap#1{\vbox to \z@{\vss\hbox to \z@{\hss#1\hss}\vss}}

\catcode`\@=12   % That's enough pretending for one day.

\catcode`\@=11

\def\mmatrix#1{\null\,\vcenter{\normalbaselines\m@th
                       \ialign{\hfil$##$\hfil&&\ \hfil$##$\hfil\crcr
                       \mathstrut\crcr\noalign{\kern-\baselineskip }
                       #1\crcr\mathstrut\crcr\noalign{\kern-\baselineskip}}}\,}
\catcode`\@=12

\def\diagram{%
\def\normalbaselines{\baselineskip15pt\lineskip3pt\lineskiplimit3pt}%
                 \mmatrix
}

\def\ln#1{\llap{$\scriptstyle #1$}}  %left name
\def\rn#1{\rlap{$\scriptstyle #1$}}  %right name

%%%%%%%%%%%%%%%%%%%%%%%%%%%%%%%%%%%%%%%%%%%%%%%%%%%%%%%%%%%%%%%%

\input amssym.def
\input amssym.tex

%%%%%%%%%%%%%%%%%%%%%%%%%%%%%%%%%%%%%%%%%%%%%%%%%%%%%%%%%%%%%%%%

%%%%%%%%%%%%%%%%%%%%%%%%%%%%%%%%%%%%%%%%%%%%%%%%%%%%%% TeX code
\catcode`\@=11
%% see \rightarrowfill, \leftarrowfill,
%%     \overrightarrow, \overleftarrow
%% of TeX

\def\Rightarrowfill{ $\m@th \mathord =\mkern -6mu\cleaders
               \hbox {$\mkern -2mu\mathord =\mkern -2mu$}
                  \hfill \mkern -6mu\mathord \Rightarrow $}
\def\Leftarrowfill{$\m@th \mathord \Leftarrow \mkern -6mu\cleaders
                     \hbox {$\mkern -2mu\mathord =\mkern -2mu$}
                     \hfill \mkern -6mu\mathord =$}
\def\Leftrightarrowfill{$\m@th \mathord \Leftarrow \mkern -6mu\cleaders
                     \hbox {$\mkern -2mu\mathord =\mkern -2mu$}
                     \hfill \mkern -6mu\mathord \Rightarrow $}

\def\leftrightarrowfill{$\m@th \mathord \leftarrow \mkern -6mu\cleaders
                     \hbox {$\mkern -2mu\mathord -\mkern -2mu$}
                     \hfill \mkern -6mu\mathord \rightarrow $}

\def\underrightarrow#1{\vtop {\m@th \ialign {##\crcr
                             $\hfil \displaystyle { #1 }\hfil $\crcr
                             \noalign {\kern 1pt \nointerlineskip }
                            \rightarrowfill \crcr \noalign {\kern 0pt }}}}

\def\underleftarrow#1{\vtop {\m@th \ialign {##\crcr
                             $\hfil \displaystyle { #1 }\hfil $\crcr
                             \noalign {\kern 1pt \nointerlineskip }
                            \leftarrowfill \crcr \noalign {\kern 0pt }}}}

\def\underleftrightarrow#1{\vtop {\m@th \ialign {##\crcr
                             $\hfil \displaystyle { #1 }\hfil $\crcr
                             \noalign {\kern 0pt \nointerlineskip }
                            \leftrightarrowfill \crcr \noalign {\kern 0pt }}}}

\def\overleftrightarrow#1{\vbox {\m@th \ialign {##\crcr
                           \noalign {\kern 3pt } \leftrightarrowfill \crcr
                           \noalign {\kern 3pt \nointerlineskip }
                           $\hfil \displaystyle { #1 }\hfil$\crcr }}}

\catcode`\@=12
%%%%%%%%%%%%%%%%%%%%%%%%%%%%%%%%%%%%%%%%%%%%%%%%%%%%%%%%%%%
                                                                  %

%%%%%%%%%%%%%%%%%%%%%%%%%%%%%%%%%%%%%%%%%%%%%%%%%%%%%%%%%%%%%%%%%%%%%%%%%%
%%%%%%% Lenght Variable Arrows
%%%%%%% this for the construction

\def\rightto#1{\hbox to#1pt{\rightarrowfill}}
\def\leftto#1{\hbox to#1pt{\leftarrowfill}}
\def\leftrightto#1{\hbox to#1pt{\leftrightarrowfill}}

\def\mapstoto#1{\mapstochar\mkern -4mu\hbox to#1pt{\rightarrowfill}}
\def\rmapstoto#1{\hbox to#1pt{\leftarrowfill}\mkern -4mu\mapstochar}

\def\hookrightto#1{\lhook\mkern -8mu\hbox to#1pt{\rightarrowfill}}
\def\hookleftto#1{\hbox to#1pt{\leftarrowfill}\mkern -8mu\rhook}

\def\twoheadrightto#1{\hbox to#1pt{\rightarrowfill}\mkern -20mu\rightarrow }
\def\twoheadleftto#1{\leftarrow\mkern -20mu\hbox to#1pt{\leftarrowfill}}

\def\hooktwoheadrightto#1{\lhook\mkern -8mu%
                                       \hbox to#1pt{\rightarrowfill}%
                                       \mkern -20mu\rightarrow }
\def\hooktwoheadleftto#1{\leftarrow\mkern -20mu%
                                      \hbox to#1pt{\leftarrowfill}%
                                      \mkern -8mu\rhook}

%%%%%%% user Arrows
%%
%% these arrows again call for the lenght in pt
%% R=right L=left
%% h=hook  t=twohead   to=to

\def\R#1{\mathop{\rightto{#1}}\limits}

\def\U#1{\left\uparrow\vbox to#1pt{}\right.}
\def\D#1{\left\downarrow\vbox to#1pt{}\right.}

%%
%% here w=big-width W=Big ww=bigg WW=Bigg

\def\W{25.5}   \def\H{12.5}

\def\point{{\scriptscriptstyle\bullet}}
\def\lra{\longrightarrow}
\def\ra{\rightarrow}
\def\id{{\rm id}}
\def\Hom{{\rm Hom}}
\def\tot{{\rm tot}}
\def\r{\right}
\def\l{\left}
\def\isom{\cong}

\bigskip

\input amssym.def
\input amssym.tex
\input epsf
%%%%%%%%%%%%%%%%%%%%%%%%%%%%%%%%%%%%%%%%%%%%%%%%%%%%%%%%%%%%%%%%%%%%%%%%%%%%%%%%
%%%	DEFINIZIONE ABSTRACT
%%%%%%%%%%%%%%%%%%%%%%%%%%%%%%%%%%%%%%%%%%%%%%%%%%%%%%%%%%%%%%%%%%%%%%%%%%%%%%%%%

\def\abstra{ \ifnum\pageno<0\else\pageno=-1\fi
          \vfill\eject
          \ifodd\pageno{}\else\quad\eject\fi
          \global\inichapter=\pageno
          \xdef\chaptertitle{Abstract}
          \xdef\sectiontitle{Abstract}
          \line{\hfill\chaptertitlefont Abstract}
          \vskip 30pt
           }
%%%%%%%%%%%%%%%%%%%%%%%%%%%%%%%%%%%%%%%%%%%%%%%%%%%%%%%%%%%%%%%%%%%%%%%%%%%%%%%%%

%% per la data di sistema

\def\oggi{\number\day\space\ifcase\month\or
       gennaio\or febbraio\or marzo\or aprile\or maggio\or giugno\or
       luglio\or agosto\or settembre\or ottobre\or novembre\or dicembre\fi
        \ \number\year}

\def\today{\ifcase\month\or
               January\or February\or March\or April\or May\or June\or
               July\or August\or September\or October\or November\or December\fi
               \space\number\day, \number\year}

\def\aujourdhui{\number\day\space\ifcase\month\or
               Janvier\or F\'evrier\or Mars\or Avril\or May\or Juin\or
               Juillet\or Ao\^ut\or Septembre\or Octobre\or Novembre\or 
               D\'ecembre\fi
               \ \number\year}

%%%%%%%%%%%%%%%%%%%%%%%%%%%%%%%%%%%%%%%%%%%%%%%%%%%%%%%%%%%%%%%%%%%%%%%%%%

%%%%%%%%%%%%%%%%%%%%%%%%%%%%%%%%%%%%%%%%%%%%%%%%%%%%%%%%%%%%%%%%%%
%% diagrammi 
%%%%%%%%%%%%%%%%%%%%%%%%%%%%%%%%%%%%%%%%%%%%%%%%%%%%%%%%%%%%%%%%%%
%%% File: Laps.TeX
\catcode`\@=11   % Let's pretend @ is a letter
%         Vertical `laps'; cf. \llap and \rlap
\long\def\ulap#1{\vbox to \z@{\vss#1}}
\long\def\dlap#1{\vbox to \z@{#1\vss}}

%         And centered horizontal and vertical `laps'
\def\xlap#1{\hbox to \z@{\hss#1\hss}}
\long\def\ylap#1{\vbox to \z@{\vss#1\vss}}

%         And a `lap' centered on its midpoint

%         And the not-long definitions
\def\ulap#1{\vbox to \z@{\vss#1}}
\def\dlap#1{\vbox to \z@{#1\vss}}
\def\hlap#1{\hbox to \z@{\hss#1\hss}}
\def\vlap#1{\vbox to \z@{\vss\hbox{#1}\vss}}
\def\clap#1{\vbox to \z@{\vss\hbox to \z@{\hss#1\hss}\vss}}

\catcode`\@=12   % That's enough pretending for one day.

\catcode`\@=11

\def\mmatrix#1{\null\,\vcenter{\normalbaselines\m@th
                    \ialign{\hfil$##$\hfil&&\ \hfil$##$\hfil\crcr
                    \mathstrut\crcr\noalign{\kern-\baselineskip }
                    #1\crcr\mathstrut\crcr\noalign{\kern-\baselineskip}}}\,}

\catcode`\@=12

\def\diagram{%
         \def\normalbaselines{\baselineskip15pt\lineskip3pt\lineskiplimit3pt}%
         \mmatrix
}

%%%%%%%%%%%%%%%%%%%%%%%%%%%%%%%%%%%%%%%%%%%%%%%%%%%%%%%%%%%%%%%%%%
%% frecce 
%%%%%%%%%%%%%%%%%%%%%%%%%%%%%%%%%%%%%%%%%%%%%%%%%%%%%% TeX code
\catcode`\@=11
%% see \rightarrowfill, \leftarrowfill,
%%     \overrightarrow, \overleftarrow
%% of TeX

\def\Rightarrowfill{ $\m@th \mathord =\mkern -6mu\cleaders
       \hbox {$\mkern -2mu\mathord =\mkern -2mu$}
          \hfill \mkern -6mu\mathord \Rightarrow $}
\def\Leftarrowfill{$\m@th \mathord \Leftarrow \mkern -6mu\cleaders
             \hbox {$\mkern -2mu\mathord =\mkern -2mu$}
             \hfill \mkern -6mu\mathord =$}
\def\Leftrightarrowfill{$\m@th \mathord \Leftarrow \mkern -6mu\cleaders
             \hbox {$\mkern -2mu\mathord =\mkern -2mu$}
             \hfill \mkern -6mu\mathord \Rightarrow $}

\def\leftrightarrowfill{$\m@th \mathord \leftarrow \mkern -6mu\cleaders
             \hbox {$\mkern -2mu\mathord -\mkern -2mu$}
             \hfill \mkern -6mu\mathord \rightarrow $}

\def\underrightarrow#1{\vtop {\m@th \ialign {##\crcr
                     $\hfil \displaystyle { #1 }\hfil $\crcr
                     \noalign {\kern 1pt \nointerlineskip }
                    \rightarrowfill \crcr \noalign {\kern 0pt }}}}

\def\underleftarrow#1{\vtop {\m@th \ialign {##\crcr
                     $\hfil \displaystyle { #1 }\hfil $\crcr
                     \noalign {\kern 1pt \nointerlineskip }
                    \leftarrowfill \crcr \noalign {\kern 0pt }}}}

\def\underleftrightarrow#1{\vtop {\m@th \ialign {##\crcr
                     $\hfil \displaystyle { #1 }\hfil $\crcr
                     \noalign {\kern 0pt \nointerlineskip }
                    \leftrightarrowfill \crcr \noalign {\kern 0pt }}}}

\def\overleftrightarrow#1{\vbox {\m@th \ialign {##\crcr
                   \noalign {\kern 3pt } \leftrightarrowfill \crcr
                   \noalign {\kern 3pt \nointerlineskip }
                   $\hfil \displaystyle { #1 }\hfil$\crcr }}}

\catcode`\@=12
%%%%%%%%%%%%%%%%%%%%%%%%%%%%%%%%%%%%%%%%%%%%%%%%%%%%%%%%%%%
%%%%%%% Lenght Variable Arrows
%%%%%%% this for the construction

\def\rightto#1{\hbox to#1pt{\rightarrowfill}}
\def\leftto#1{\hbox to#1pt{\leftarrowfill}}
\def\leftrightto#1{\hbox to#1pt{\leftrightarrowfill}}

\def\mapstoto#1{\mapstochar\mkern -4mu\hbox to#1pt{\rightarrowfill}}
\def\rmapstoto#1{\hbox to#1pt{\leftarrowfill}\mkern -4mu\mapstochar}

\def\hookrightto#1{\lhook\mkern -8mu\hbox to#1pt{\rightarrowfill}}
\def\hookleftto#1{\hbox to#1pt{\leftarrowfill}\mkern -8mu\rhook}

\def\twoheadrightto#1{\hbox to#1pt{\rightarrowfill}\mkern -20mu\rightarrow }
\def\twoheadleftto#1{\leftarrow\mkern -20mu\hbox to#1pt{\leftarrowfill}}

\def\hooktwoheadrightto#1{\lhook\mkern -8mu%
                               \hbox to#1pt{\rightarrowfill}%
                               \mkern -20mu\rightarrow }
\def\hooktwoheadleftto#1{\leftarrow\mkern -20mu%
                              \hbox to#1pt{\leftarrowfill}%
                              \mkern -8mu\rhook}

%%%%%%% user Arrows
%%
%% these arrows again call for the lenght in pt
%% R=right L=left
%% h=hook  t=twohead   to=to

\def\R#1{\mathop{\rightto{#1}}\limits}

\def\U#1{\left\uparrow\vbox to#1pt{}\right.}
\def\D#1{\left\downarrow\vbox to#1pt{}\right.}

%%
%% here w=big-width W=Big ww=bigg WW=Bigg

\def\W{25.5}   \def\H{12.5}

%% to use as in $\R\W_{}^{}$

%%%%%% we wont also define Up and Down h/t/ht/to arrows
%% but giving the right dimensions, easy for ortho-transf, so
%% read part of rotate.tex
\newdimen\rotdimen
\def\vspec#1{\special{ps:#1}}%  passes #1 verbatim to the output
\def\rotstart#1{\vspec{gsave currentpoint currentpoint translate
        #1 neg exch neg exch translate}}
	   % #1 can be any origin-fixing transformation
\def\rotfinish{\vspec{currentpoint grestore moveto}}% gets back in synch

\def\rotr#1{\rotdimen=\ht#1\advance\rotdimen by\dp#1%
        \hbox to\rotdimen{\hskip\ht#1%
        \vbox to\wd#1{\rotstart{90 rotate}%
        \box#1\vss}\hss}\rotfinish}%

\def\rotr#1{\rotdimen=\ht#1\advance\rotdimen by\dp#1%
        \hbox to\rotdimen{\hskip\ht#1
        $\vcenter to\wd#1{\rotstart{90 rotate}%
        \box#1\vss}$\hss}\rotfinish}%

\def\rotl#1{\rotdimen=\ht#1\advance\rotdimen by\dp#1%
        \hbox to\rotdimen{\vbox to\wd#1{\vskip\wd#1\rotstart{270 rotate}%
        \box#1\vss}\hss}\rotfinish}%

%% and make the definitions
\newbox\horarr
\def\Uh#1{\setbox\horarr=\hbox{$\hookleftto{#1}$}\rotr{\horarr}}
\def\Dh#1{\setbox\horarr=\hbox{$\hookrightto{#1}$}\rotr{\horarr}}
\def\Ut#1{\setbox\horarr=\hbox{$\twoheadrightto{#1}$}\rotl{\horarr}}
\def\Dt#1{\setbox\horarr=\hbox{$\twoheadleftto{#1}$}\rotl{\horarr}}
\def\Uht#1{\setbox\horarr=\hbox{$\hooktwoheadleftto{#1}$}\rotr{\horarr}}
\def\Dht#1{\setbox\horarr=\hbox{$\hooktwoheadrightto{#1}$}\rotr{\horarr}}
\def\Uto#1{\setbox\horarr=\hbox{$\mapstpto{#1}$}\rotl{\horarr}}
\def\Dto#1{\setbox\horarr=\hbox{$\rmapstoto{#1}$}\rotl{\horarr}}

%%%%%%%%%%%%%%%%%%%%%%%%%%%%%%%%%%%%%%%%%%%%%%%%%%%%%%%%%%%%%%%%%%%%%
%% PSx\TeX utils by Maurizio Cailotto (maurizio@math.unipd.it)
%%
%% centred rotations of #1 degree of #2
%% end
%% centred #1 x axis, #2 y axis deformation of #3
%% but
%% the results are of zero dimensions, for the moment!

%% Centred ROTations
%% this is the skeleton

%% Centred DEFormations
%% this is the skeleton
\def\cdef#1#2#3{
        \vbox to0pt {\vss \hbox to0pt {\hss #3 \hss} \vss}
}
%% this is an experimental version of the previous skeleton
%% taking care of dimensions (here it's easy)
\newbox\deform
\newdimen\altezza
\newdimen\profondita
\newdimen\larghezza
\def\cdef#1#2#3{
        \setbox\deform=\hbox{#3}
        \altezza=\ht\deform%
        \profondita=\dp\deform%
        \larghezza=\wd\deform%
       \ifnum #1>0 \multiply\larghezza by #1 \else\multiply\larghezza by -#1\fi%
       \ifnum #2>0 \multiply\altezza by #2 \else\multiply\altezza by -#2\fi%
       \ifnum #2>0 \multiply\profondita by #2 \else\multiply\profondita
by -#2\fi%
       \advance\altezza by \profondita%
        \hbox to\larghezza{\hss\lower\profondita
        \vbox to \altezza{\vss{
        \vbox to0pt {\vss \hbox to0pt {\hss #3 \hss} \vss}
        \vss}}
        \hss}
        }

%%%%%%%%%%%%%%%%%%%%%%%%%%%%%%%%%%%%%%%%%%%%%%%%%%%%%%%%%%%%%%%%%

%%%%%%%%%%%%%%%%%%%%%%%%%%%%%%%%%%%%%%%%%%%%%%%%%%%%%%%%%%%%%%%%%
%% matrici piccole 
\catcode`\@=11
\def\smalleqalign#1{\null \,\vcenter{\openup \jot \m@th %
       \ialign {\strut \hfil $\scriptstyle {##}$%
                &$\scriptstyle {{}##}$\hfil \crcr #1\crcr }}\,}
\def\smallmatrix#1{\null \,\vcenter {\baselineskip=6pt \m@th 
             \ialign {\hfil $\scriptstyle ##$\hfil &&\  
             \hfil $\scriptstyle ##$\hfil \crcr 
             \mathstrut \crcr \noalign {\kern -7pt } #1\crcr 
             \mathstrut \crcr \noalign {\kern -7pt }}}\,}
\def\verysmallmatrix#1{\null \,\vcenter {\baselineskip=4pt \m@th 
             \ialign {\hfil $\scriptscriptstyle ##$\hfil &&\  
             \hfil $\scriptscriptstyle ##$\hfil \crcr 
             \mathstrut \crcr \noalign {\kern -5pt } #1\crcr 
             \mathstrut \crcr \noalign {\kern -5pt }}}\,}

\catcode`\@=12

%%%%%%%%%%%%%%%%%%%%%%%%%%%%%%%%%%%%%%%%%%%%%%%%%%%%%%%%%%%%%%%%%
%% abbreviazioni 

\def\ds{\displaystyle}

\def\sss{\scriptscriptstyle}

\def\ov{\overline}

\def\r{\right}
\def\l{\left}

\def\phi{\varphi}

\def\theta{\vartheta}

\def\rho{\varrho}

\def\epsilon{\varepsilon}

\def\NoBlackBoxes{\global\overfullrule 0pt}
\NoBlackBoxes

\def\ra{\mathop{\rightarrow}\limits}

\def\lra{\mathop{\longrightarrow}\limits}

\def\ln#1{\llap{$\scriptstyle #1$}}  %left name
\def\rn#1{\rlap{$\scriptstyle #1$}}  %right name

\def\leq{\leqslant}

%% parentesi doppie

%% font matematici: bold, gotic, calligrafico

%%%%%%%%%%%%%%%%%%%%%%%%%%%%%%%%%%%%%%%%%%%%%%%%%%%%%%%%%%%%%%%%%%%%%%%%%%
% The following allows the use of Ralph Smith's Formal Script symbols
% in Plain TeX documents.  Use \scr like \cal.
% Set the font sizes and restore the `at' clauses if you want them bigger.
% You can use this method in LaTeX, but only at one basic size.
% If you need symbols in LaTeX titles, captions, etc., work it out or ask
% a LaTeXpert.

\font\tenscr=rsfs10 % scaled \magstep1
\font\sevenscr=rsfs7 % scaled \magstep1
\font\fivescr=rsfs5 % scaled \magstep1
\skewchar\tenscr='177 \skewchar\sevenscr='177 \skewchar\fivescr='177
\newfam\scrfam \textfont\scrfam=\tenscr \scriptfont\scrfam=\sevenscr
\scriptscriptfont\scrfam=\fivescr
\def\scr{\fam\scrfam}
%%%%%%%%%%%%%%%%%%%%%%%%%%%%%%%%%%%%%%%%%%%%%%%%%%%%%%%%%%%%%%%%%%%%%%%%%%%

\def\c#1{ {\cal #1} }
\def\c#1{ {\scr #1} }
%%%%%%%%%%%%%%%%%%%%%%%%%%%%%%%%%%%%%%%%%%%%%%%%%%%%%%%%%%%%%%%%%%%%%%%%%%%

%%%%%%%%%%%%%%%%%%%%%%%%%%%%%%%%%%%%%%%%%%%%%%%%%%%%%%%%%%%%%%%%%%%%%%%%%%%
%% lettere greche (italiche o \c )
\def\cOmega{{\mit\Omega}}
\def\cTheta{{\mit\Theta}}

%%%%%%%%%%%%%%%%%%%%%%%%%%%%%%%%%%%%%%%%%%%%%%%%%%%%%%%%%%%%%%%%%%%%%%%%%%%
%% nomi 

\def\Hom{{\rm Hom}}

\def\cHom{{{\c H}{\mit o m}}}

\def\O{{\rm O}}

\def\id{{\rm id}}

\def\bydef{\mathrel{:=}}

\def\isom{\cong}

%%%%%%%%%%%%%%%%%%%%%%%%%%%%%%%%%%%%%%%%%%%%%%%%%%%%%%%%%%%%%%%%%%%%%%%%%%%
%% Limiti
%% vedere comunque quelli AMS
%% la prima serie richiede le frecce estendibili

%% e derivati

%%%%%%%%%%%%%%%%%%%%%%%%%%%%%%%%%%%%%%%%%%%%%%%%%%%%%%%%%%%%%%%%%%%%%%%%%%%
%%	AGGIUNTE DI LUISA AL FILE DI MAURIZIO
%%%%%%%%%%%%%%%%%%%%%%%%%%%%%%%%%%%%%%%%%%%%%%%%%%%%%%%%%%%%%%%%%%%%%%%%%%%

%%%%%%%%%%%%%%%%%%%%%%%%%%%%%%%%%%%%%%%%%%%%%%%%%%%%%%%%%%%%%%%%%%%%%%%%%%%
%%%%%%%%%%%%%%%%%%%%%%%%%%%%%%%%%%%%%%%%%%%%%%%%%%%%%%%%%%%%%%%%%%%%%%%%%%%

 %%%%%% SERVE PER SALVARE LA ??? %%%%

%%%%%%%%%%%%%%%%%%%%%%%%%%%%%%%%%%%%%%%%%%%%%%% symbols

\def\id{{\rm id}}
\def\bydef{\mathrel{:=}}

\def\Diff{{\rm Diff}}

\def\epsilon{\varepsilon}

\def\point{{\scriptscriptstyle\bullet}}
\def\DR{{\rm DR}}
\def\tDR{\widetilde{\rm DR}}

\def\bR{{\bf R}}
\def\bL{{\bf L}}

\def\gr{{\rm gr}}
\def\Gr{{\rm Gr}}

\def\tot{{\rm tot}}

\def\circdot{{\circ\kern-3.5pt\lower .8pt\hbox{$\cdot$}}}

%%%%%%%%%%%%%%%%%%%%%%%%%%%%%%%%%%%%%%%%%%%%%%%%%%%%%%%%%%

\abovedisplayskip=6pt plus2pt minus 4pt
\belowdisplayskip=6pt plus2pt minus 4pt
\abovedisplayshortskip=0pt plus2pt
\belowdisplayshortskip=2pt plus1pt minus 1pt

%%%%%%%%%%%%%%%%%%%%%%%%%%%%%%%%%%%%%%%%%%%%%%%%%%%%%%%%%%%%%%%%%%%%%%%%%%%%%%%%%%%%%%%%%%%%%%%%%%%%%%%%%%%%%%%%%%%%%%%%%%%%%%%%%%%%%%%%%%%%%%%%%%%%%%%%%%%%%%%%%%%%%%%%%%%%%%%%%%%%%%%%%%%%%%%%%%%%%%%%%%%%%%%%%%%%%%%%%%%%%%%%%%%%%%%%%%%%%%%%%%%%%%%%%%%%%%%%%%%%%%%%%%%%%%%%%%%%%%%%%%%%%%%%%5

\def\I(#1){{\rm Ind(\c #1)}}
\def\P(#1){{ \rm Pro(\c #1)}}
\def\bs{{\bigskip}}
\def\O(#1){{\rm {Ob}(\c #1)}}
\def\o(#1){{\rm {Ob}(\rm #1)}}

\def\coh{\rm {coh}}

\def\qcoh{\rm {qcoh}}

\def\opp{\circ}

\def\sr{{\sss r}}

\input xy
\xyoption{all}
%\input bull.sty
%%%%%%%%%%%%%%%%%%%%%%%%%%%%%%%%%%%%%%%%%%%%%%
%%%%%%%%%		TITOLO				%%%%%%%%%%%%%%%%%%%%%
%%%%%%%%%%%%%%%%%%%%%%%%%%%%%%%%%%%%%%%%%%%%%%
\def\nomefile{art}
%%%%%%%%%%%%%%%%%%%%%%%%%%%%%%%%%%%%%%%%%%%%%%
\title{{On derived categories of differential complexes}}
%%%%%%%%%%%%%%%%%%%%%%%%%%%%%%%%%%%%%%%%%%%%%%

%%%%%%%%%%%%%%%%%%%%%%%%%%%%%%%%%%%%%%%%%%%%%%
\author{Luisa Fiorot
\footnote{}{\rm 
2000 AMS classification 14F30. 
Partially supported by PGR "CPDG021784" (University of Padova).}
}
%%%%%%%%%%%%%%%%%%%%%%%%%%%%%%%%%%%%%%%%%%%%%%

%%%%%%%%%%%%%%%%%%%%%%%%%%%%%%%%%%%%%%%%%%%%%%
%%%%%%%%%		ABSTRACT			%%%%%%%%%%%%%%%%%%%%
%%%%%%%%%%%%%%%%%%%%%%%%%%%%%%%%%%%%%%%%%%%%%%
\resume{
Cet article est  consacr\'e \`a la comparaison 
entre diff\'erentes cat\'egories localis\'ees de complexes
diff\'erentiels.
%\endgraf
Nous prouvons que le functeur canonique de
la cat\'egorie des complexes diff\'erentiels d'ordre  un
(d\'efinie par  Herrera et Lieberman) \`a valeurs dans la cat\'egorie des 
complexes diff\'erentiels (d'ordre arbitraire, definie par M. Saito),
localis\'ees par rapport \`a une bonne notion de quasi-isomorphismes, 
est une \'equivalence de cat\'egories.
%\endgraf
En suite nous prouvons un r\'esultat analogue pour une version filtr\'ee des 
cat\'egories pr\'ec\'edentes (d\'efinies respectivement par Du Bois et M. Saito),
localis\'ees par les quasi-isomorphismes gradu\'es.
Cet r\'esultat reponde \`a une question pos\'ee par M. Saito.
}
\abstract { 
This paper is devoted to the comparison of different localized categories
of differential complexes.
%\endgraf
The first result is that the canonical functor from the category of
complexes of differential operators of order one (defined by Herrera and Lieberman)  to
the category of differential complexes (of any order, defined by M. Saito),
both localized with respect to a suitable notion of  quasi-isomorphism, is an
equivalence of categories.
%\endgraf
Then we prove a similar result for a filtered version of the previous
categories (defined respectively by Du Bois and M.Saito), localized with
respect to graded-quasi-isomorphisms, thus answering a question posed
by M. Saito.
		}
\apririferimenti

%%%%%%%%%%%%%%%%%%%%%%%%%%%%%%%%%%%%%%%%%%%%
%%%%%%%%%		INRODUCTION			%%%%%%%%%%%%%%%%
%%%%%%%%%%%%%%%%%%%%%%%%%%%%%%%%%%%%%%%%%%%%

\intro{
The category of differential complexes  appears naturally as a ``good'' category for the
role of image of the classical De Rham functor.
We are interested in finding a purely algebraic definition for the image category
of the De Rham functor for differential modules which will permit us to develop the formalism
of the six Grothendieck operations.
\endgraf
Such a category was first introduced by Herrera-Lieberman in
their article in
{\it Inventiones Math.} of 1971.
In that paper they proposed the study of a category $C_1(\c O_X,\Diff_X)$  of complexes of
$\c O_X$-Modules with differential operators of order one (where $X$ is
a smooth algebraic or analytic variety over a field $K$ of
characteristic zero).
They also interpreted  $C_1(\c O_X,\Diff_X)$  as a category of graded modules
over a suitable graded ring $\c C^{\point}_X$ containing 
$\cOmega^{\point}_X$ as a sub-ring.
Using this interpretation they defined the functors $-\otimes_{\cOmega^\point}-$ , 
$\cHom_{\cOmega^\point}(-,-)$, $f^\ast$ and $f_\ast$.
Then they defined hyperext functors using suitable injective resolutions
and  proved a duality theorem in the proper smooth case.
\endgraf
They did not propose in that paper to localize 
$C_1(\c O_X, \Diff_X)$ with respect
to a multiplicative system as is done in the study of derived categories, although they did 
introduce a notion of  homotopy.
\endgraf
The difficulty  in the localization procedure was first pointed out by 
P. Berthelot in his book of 1974 [B] , where he showed that objects of 
$C_{1}(\c O_X,\Diff_X)$ which are quasi-isomorphic as complexes of 
abelian sheaves, may lead to non isomorphic hyperext functors.
On the other hand, we are forced to localize $C_{1}(\c O_X,\Diff_X)$, if we wish 
to obtain a triangulated category where a De Rham functor $\DR$, with 
source some derived category of ${\c D}_X$-Modules, can assume its values.
By Berthelot's remark we know that the multiplicative system of abelian
quasi-isomorphisms is not a good choice.
\endgraf
Different localizations were proposed by Philippe Du Bois who, in [DB.1], introduced
filtrations and so obtained the category $DF_1(\c O_X,\Diff_X)$, and by Morihiko Saito who,  in [S.1] and [S.2],
defined a new category of complexes $C(\c O_X,\Diff_X)$ with
differential operators (of any order) which he
localized with respect to $\tDR^{-1}_X$-quasi-isomorphism (obtaining $D(\c O_X,\Diff_X)$)
or with respect to filtered quasi-isomorphism (obtaining $DF(\c O_X,\Diff_X)$).
\endgraf
 Saito's category   $C(\c O_X,\Diff_X)$ seems to be the best choice because it is equivalent to the category
$D(\c D_X)$ via the functors $\tDR_X$ and $\tDR_X^{-1}$.
The problem is that in Saito's category an explicit  
formalism of Grothendieck operations is only partially realized; 
in fact, for example, there is no $f^\ast$ functor or internal tensor
product.
From the Herrera-Lieberman point of view, considering the category $D_1(\c O_X,\Diff_X)$ 
obtained by localizing $C_1(\c O_X,\Diff_X)$ with respect to $\tDR^{-1}_X$-quasi-isomorphism,
we obtain a category wherein the De Rham functor takes its image and where 
the Grothendieck operations are easier and more complete than those in Saito's
category $D(\c O_X,\Diff_X)$.

This work is  devoted to the comparisons between the categories
of differential complexes $D_1(\c O_X,\Diff_X)$ and $D(\c O_X,\Diff_X)$.
In the first section we recall some general definitions we need in this paper
and the notation we will use.
Then in Section 2  we develop a general result about morphisms to
a total complex in a general category of complexes.
In Section 3 we compare  Saito category $D(\c O_X,\Diff_X)$ 
with $D_1(\c O_X,\Diff_X)$.
In fact we prove that the canonical functor between the localized categories
$i_{HL,S}:D_1(\c O_X,\Diff_X)\lra D(\c O_X,\Diff_X)$ 
is an equivalence of categories.
In the last section we extend this comparison result to the filtered case proving that
the filtered Du Bois category $DF_1(\c O_X,\Diff_X)$ and Saito's 
$DF(\c O_X,\Diff_X)$ are equivalent, thus answering a question posed by Saito in
[S.1, 2.2.11].
\smallskip 
I would like to thank Prof. Francesco Baldassarri 
for having introduced me to this matter.
It is a pleasure to thank Maurizio Cailotto and Morihiko Saito for the improvements 
and suggestions they gave me in the redaction of this work.
}
\bs
%%%%%%%%%%%%%%%%%%%%%%%%%%%%%%%%%%%%%%%%%%%%%%%%%%
%%%%%%%%%%%%%%%%%%%%%%%%%%%%%%%%%%%%%%%%%%%%%%%%%%

%%%%%%%%%%%%%%%%%%%%%%%%%%%%%%%%%%%%%%%%%%%%%%%%%%
\section{Notation and definitions.}
%%%%%%%%%%%%%%%%%%%%%%%%%%%%%%%%%%%%%%%%%%%%%%%%%%

Let $X$ be a smooth separated scheme of finite type  over a field $K$ of 
characteristic zero, or a smooth analytic variety.
\smallskip

%%%%%%%%%%%%%%%%%%%%%%%%%%%%%%%%%%%%%%%%%%%%%%%%%
%% Definizione di Herrera-Liebermann differential complexes%%%%%%%%%%%%%%%%%%
%%%%%%%%%%%%%%%%%%%%%%%%%%%%%%%%%%%%%%%%%%%%%%%%%

\defi{NN}
{\rif{HLCX}{\bf Herrera-Liebermann differential complexes.}\endgraf
As in [HL, \S 2] or [B, II.5], the category $C_{1}(\c O_X,\Diff_X)$
is defined as the category of complexes of differential operators of order 
at most one,
that is:
\endgraf
	\item{i)} the objects of $C_{1}(\c O_X,\Diff_X)$ are complexes whose terms are
		$\c O_{X}$-Modules and whose differentials are differential operators
		of order less than or equal to one;
\endgraf
	\item{ii)} morphisms between such complexes are morphisms of complexes
		which are $\c O_{X}$-linear maps.
\endgraf\noindent
We denote by $C_{1}^b(\c O_X,\Diff_X)$ the full subcategory of $C_{1}(\c O_X,\Diff_X)$
whose objects are bounded complexes.
}
\smallskip

%%%%%%%%%%%%%%%%%%%%%%%%%%%%%%%%%%%%%%%%%%%%%%%%%%%

\prop{NN}
{
The category $C_{1}(\c O_X,\Diff_X)$ is equivalent to the category of graded
left $\c C^{\point}_{X}$-Modules where
$\c C^{\point}_{X}\isom
{\cOmega}_{X}^{\point{-}1}D\oplus{\cOmega}_{X}^\point$
is the ``mapping cylinder''
of the identity map of ${\cOmega}_{X}^\point$.
It is a graded $\c O_{X}$-Algebra, whose product is defined
using the wedge product of ${\cOmega}_{X}^\point$ and $D^{2}=0$,
while the structure of complex is defined by
$D\alpha=(d\alpha_{1}{+}(-1)^{i}\alpha_{2})D+d\alpha_{2}$
if $\alpha=\alpha_{1}D+\alpha_{2}$ with
$\alpha_{1}\in\cOmega_{X}^{i{-}1}$ and
$\alpha_{2}\in\cOmega_{X}^{i}$.
Therefore, the category $C_{1}(\c O_X,\Diff_X)$ has enough injectives
[HL, \S 2].
\rif{HLCXM}}
\smallskip

%%%%%%%%%%%%%%%%%%%%%%%%%%%%%%%%%%%%%%%%%%%%%%%%%
%% Definizione di omotopia 			%%%%%%%%%%%%%%%%%%%%%%%%%%%%%
%%%%%%%%%%%%%%%%%%%%%%%%%%%%%%%%%%%%%%%%%%%%%%%%%

\defi{NN}
{
A homotopy between two morphisms in $C_{1}(\c O_X,\Diff_X)$ is a
homotopy in the sense of the category of complexes of abelian
sheaves, except that the homotopy operator (of degree $-1$)  is
taken to be $\c O_{X}$-linear (see [HL, \S 2]).
}
\smallskip

%%%%%%%%%%%%%%%%%%%%%%%%%%%%%%%%%%%%%%%%%%%%%%%%%
%% Definizione di operatore differenziale		        		%%%%%%%%%%%%%%%%%%
%%%%%%%%%%%%%%%%%%%%%%%%%%%%%%%%%%%%%%%%%%%%%%%%%
\defi{NN}
{
Let  $\c D_X$ be the sheaf of differential operators on $X$ (see [Bo] for the 
definition) and $p:\c D_X\lra \c O_X$ 
the map evaluating a differential operator in 1.
In [S.2] Saito defines a differential operator $d:\c F\lra \c G$, between
two $\c O_X$-Modules $\c F$ and $\c G$, as a morphism which can be factorized 
(in a unique way) as
$$
\xymatrix{
\c F	\ar[r]^{d} \ar[rd]_{\ov d}	&
\c G					\cr
					&
\c G\otimes_{\c O_X}\c D_X \ar[u]_{id_{\c G}\otimes_{\c O_X}p}
}
$$
where $\ov d$ is an $\c O_X$-linear map and 
$\c G\otimes_{\c O_X}\c D_{X}$
is an $\c O_X$-Module for the right multiplication of $\c D_{X}$.
\endgraf
We note that the morphism  $\ov d$ induces  by extension of scalars
a morphism of $\c D_X$-Modules
$\ov d': \c F\otimes_{\c O_X}\c D_{X}
\lra \c G\otimes_{\c O_X}\c D_{X}$.
\endgraf
We use the notation $M(\c O_X,\Diff_X)$ for the additive category whose objects are $\c O_X$-Modules
and morphisms are differential operators between them.
}
\smallskip

%%%%%%%%%%%%%%%%%%%%%%%%%%%%%%%%%%%%%%
%% Definizione di DR^-1					 	%%%%%%%%%%%%%
%%%%%%%%%%%%%%%%%%%%%%%%%%%%%%%%%%%%%%
\defi{NN}
{{\bf Saito differential complexes.}\endgraf
In [S.2] Saito defines the equivalence of categories 
$$\matrix{
\tDR^{-1}_X: & M(\c O_X,\Diff_X) & \lra & M_i(\c D_X)^\sr \cr
	  & \c F	      & \longmapsto & \c F\otimes_{\c O_X}\c D_X\cr
	  & d                 & \longmapsto & \ov d'			\cr
}
$$
where $M_i(\c D_X)^\sr$ is the full subcategory of right
$\c D_X$-Modules whose objects are 
induced modules (i.e. they are of the form $\c F\otimes_{\c O_X}\c D_X$ for an
$\c O_X$-Module $\c F$). 
\endgraf
Let $C(\c O_X,\Diff_X)$ 
be the category of complexes in $M(\c O_X,\Diff_X)$.
Then the $\tDR_X^{-1}$ functor extends to a functor
$\tDR_X^{-1}:C(\c O_X,\Diff_X)\lra C(\c D_X)^\sr$. 
}
\smallskip

%%%%%%%%%%%%%%%%%%%%%%%%%%%%%%%%%%%%%%%%%%%%

\rema{NN}
{
In Definition 1.1 we have introduced the Herrera-Lieberman category $C_1(\c O_X,\Diff_X)$.
The main difference between $C_1(\c O_X,\Diff_X)$ and $C(\c O_X,\Diff_X)$ is that
Herrera-Liebermann allow only differential operators of order one and morphisms 
between complexes are $\c O_X$-linear, while Saito considers complexes with 
arbitrary differential operators and morphisms 
between complexes given by  differential operators.
We observe that the map $p:\c D_X\lra \c O_X$ is a differential operator but 
it is not of finite order.
\endgraf
There is a natural functor
$\lambda_1 :C_1(\c O_X,\Diff_X)\lra C(\c O_X,\Diff_X)$ 
which sends objects of the first category into themselves regarded as 
objects of $C(\c O_X,\Diff_X)$.
This functor is not faithful.
}
\smallskip

%%%%%%%%%%%%%%%%%%%%%%%%%%%%%%%%%%%%%%%%%%%%

\defi{NN}
{
We define  $D_1(\c O_X,\Diff_X)$ to be the category obtained localizing
the category $C_1(\c O_X,\Diff_X)$ of Herrera and Lieberman with respect to
the multiplicative system of morphisms 
$$
S_{1,{\tDR}^{-1}_X}
:=\{f\in C_1(\c O_X,\Diff_X)| \ {\tDR}^{-1}_{X}\opp i_{HL,S}(f)
\hbox { \ is\ a\ quasi-isomorphism\ in\ }
D(\c D_{X})^\sr \}
$$
called the system of ${\tDR_X}^{-1}$-quasi-isomorphisms.
We refer to this category as the Herrera-Lieberman localized category.
}
\smallskip

%%%%%%%%%%%%%%%%%%%%%%%%%%%%%%%%%%%%%%%%%%%%%%%%%
%% Definizione di Saito differential complexes		%%%%%%%%%%%%%%%%%%%%%
%%%%%%%%%%%%%%%%%%%%%%%%%%%%%%%%%%%%%%%%%%%%%%%%%

\defi{NN}
{
We define the category 
$D(\c O_X,\Diff_X)$ 
by localizing  the category 
$C(\c O_X,\Diff_X)$ 
with respect to the multiplicative system
$$
S_{{\tDR}^{-1}_X}
:=\{f\in C(\c O_X,\Diff_X)| \ {\tDR}^{-1}_{X}(f)
\hbox { \ is\ a\ quasi-isomorphism\ in\ }
D(\c D_{X})^\sr \}.
$$
}
\smallskip

%%%%%%%%%%%%%%%%%%%%%%%%%%%%%%%%%%%%%%%%%%%%%%%%%

\rema{NN}
{ 
The natural functor $\lambda_1$ respects
${\tDR}^{-1}_{X}$-quasi-isomorphisms
so it defines a functor (which we again denote by $\lambda_1$)
$$
\lambda_1: D_1(\c O_X,\Diff_X)\lra D(\c O_X,\Diff_X).
$$
It seems to be not straightforward to prove that this functor is
fully faithful.
\endgraf
\smallskip
We recall that Saito proved in [S.2] that the usual De Rham functor 
(for right $\c D_X$-Modules)
$$
\DR_{X}
:={-\otimes^{\bL}_{\c D_{X}}{\c O_X}}:
D^b(\c D_{X})^\sr \lra  D^b(K_X)
$$
factors as 
$$
\xymatrix{
D^b(\c D_{X})^\sr\ar[rd]_{\DR_{X}}\ar[r]^{\widetilde{\DR}_{X}}
& D^b(\c O_X, \Diff_X)\ar[d]	\\
& D^b(K_X)\\
}
$$ 
where the vertical arrow is the functor obtained by forgetting all 
structure of an object in
$D^b(\c O_X, \Diff_X)$ 
but that of 
$K_X$-Module;
 and $\widetilde{\DR}_{X}$ is the functor
$$
\matrix{
\widetilde{\DR}_{X}:& 
D^b(\c D_{X})^{\sr}  &  \lra & D^b(\c O_X, \Diff_X)\cr
&\c M^{\point}&\longmapsto& \c M^{\point}\otimes_{\c
O_X}^{\point}{\cTheta^{\point}_{X}}.\cr}
\leqno{\bf (\NNN)}\rif{DR}
$$
($\cTheta^i_X=\wedge^{-i}\cTheta_X$)
which could be extended to unbounded complexes.
\endgraf
For 
$\c M^{\point}\in D^b(\c D_{X})^{\sr}$ 
we have the following
three descriptions of $\DR_{X}(\c M^{\point})$:
$$
\eqalign{
\DR_{X}(\c M^{\point}) 
&  =\c M^{\point}\otimes^{\bL}_{\c D_{X}}{\c O_X}\cr
&  =\c M^{\point}\otimes^{\point}_{\c O_X}\cTheta^{\point}_{X}\cr
&  =\bR\cHom_{\c D_{X}}(\omega_{X},\c M^{\point})[n]\cr }
$$ 
in $D^b(K_X)$ where $n= \dim X$.
\endgraf
We observe that the functor $\tDR_X$, extended to unbounded complexes, 
also factors  through $D_1(\c O_X,\Diff_X)$, so we obtain the commutative
diagram
$$
\xymatrix{
D(\c D_{X})^\sr   \ar[r]^{{\tDR}_{1,X}}   
\ar[rd]_{{\tDR}_{X}} &
D_1(\c O_X,\Diff_X)   \ar[d]^{\lambda_1} \\
&  D(\c O_X,\Diff_X).\\
}
$$
The composition
${\tDR}_{1,X}{\tDR}_{X}^{-1}$ 
defines a functor
$D(\c O_X,\Diff_X)\lra D_1(\c O_X,\Diff_X)$.
In the sequel, we will  prove that $\lambda_1$ and  
${\tDR}_{1,X}{\tDR}_{X}^{-1}$ are quasi-inverses of each other and so define an equivalence of categories.
}
\smallskip

\theo{NN}
{\kern -16pt{\rm (Saito[S.2])}  The functors
$\tDR_X$ and $\tDR^{-1}_X$ are equivalences between the categories
$D(\c D_X)^\sr$ and $D(\c O_X,\Diff_X)$.
}
\smallskip

\eject
%%%%%%%%%%%%%%%%%%%%%%%%%%%%%%%%%%%%%%%%%%%%%%%%%%%%%%%%%%%%%%%%%%%%%%%%%%%%%
\section{Morphisms between a complex and a total complex.}
%%%%%%%%%%%%%%%%%%%%%%%%%%%%%%%%%%%%%%%%%%%%%%%%%%%%%%%%%%%%%%%%%%%%%%%%%%%%%

{\NN} {\bf Notation.}
In this section, we consider an additive category $\c A$, the category of complexes 
and the category of na{\"\i}f bounded bicomplexes on it 
(na{\"\i}f means that the differentials commute, 
and bounded means that in any anti-diagonal only a finite number of 
terms are not isomorphic to zero). 
We want to describe the morphisms of complexes between a complex 
and the total complex associated to a bicomplex. 
\smallskip

\lemm{NN} 
{
Let $I^{\point,\point}$ be a bicomplex with commuting differentials 
$d^{\prime}_{I}$ and $d^{\prime\prime}_{I}$. 
We define the bicomplex $\widetilde{I}^{\point,\point}$ 
in the following way:
for any $p$ and $q$ let
$\widetilde{I}^{p,q}={I}^{p,q}$,
$d^{\prime\,p,q}_{\widetilde{I}}=
(-1)^{q}d^{\prime\,p,q}_{I}$
and
$d^{\prime\prime\,p,q}_{\widetilde{I}}=
(-1)^{p}d^{\prime\prime\,p,q}_{I}$. 
Then: 
\item{$(i)$} 
the functor $\;\widetilde{}\;$ sending $I^{\point,\point}$
to $\widetilde{I}^{\point,\point}$ is an automorphism of the
category of bicomplexes, and $\skew4\widetilde{\widetilde{I}}=I$ 
for any bicomplex $I$.
\item{$(ii)$} 
the canonical map
$\sigma_{I}^{\point,\point}:
I^{\point,\point}\lra\widetilde{I}^{\point,\point}$
defined by
$\sigma_{I}^{p,q}=(-1)^{pq}\id_{I^{p,q}}$
is an isomorphism of bicomplexes
and defines an isomorphism of functors 
$\sigma:\id\ra\;\widetilde{}\;$. 
} 
\smallskip

%%%%%%%%%%%%%%%%%%%%%%%%%%%%%%% ~ Hom^{\point}

\coro{NN} 
{
Let $A^{\point}$, resp;
$B^{\point,\point}$, be a complex resp. a bicomplex in $\c A$.
Then we have a commutative diagram of canonical isomorphisms
of bicomplexes 
$$ 
\diagram{
\Hom_{}^{s}(A^{\point},B^{\point,q}) 
&\R\W^{\sigma^{s,q}_{B\ast}} 
&\Hom_{}^{s}(A^{\point},\widetilde{B}^{\point,q}) \cr 
\ln{\sigma^{s,q}_{\Hom}}\D\H 
&&\D\H\rn{\widetilde\sigma^{s,q}_{\Hom}} \cr 
\widetilde\Hom_{}^{s}(A^{\point},B^{\point,q})
&\R\W_{\widetilde\sigma^{s,q}_{B\ast}} 
&\widetilde\Hom_{}^{s}(A^{\point},\widetilde{B}^{\point,q})
}$$ 
where 
$$
\widetilde\Hom_{}^{s}
(A^{\point},{B}^{\point,q})_{s,q}=
\left(\Hom_{}^{s}
(A^{\point},{B}^{\point,q})_{s,q}\right){\widetilde{\ }}
$$
and
$$
\widetilde\Hom_{}^{s}
(A^{\point},\widetilde{B}^{\point,q})_{s,q}=
\left(\Hom_{}^{s}
(A^{\point},\widetilde{B}^{\point,q})_{s,q}\right){\widetilde{\ }} \; ,
$$
the horizontal morphisms are 
$$ 
\sigma_{A\ast}=
\Hom_{}^{\point}(\id_{A^{\point}},\sigma_{B})
\quad\hbox{and}\quad
\widetilde\sigma_{B\ast}=
\widetilde\Hom_{}^{\point}(\id_{A^{\point}},\sigma_{B})\; ,
$$ 
and the vertical ones are 
$$ 
\sigma^{s,q}_{\Hom}=\sigma_{\Hom_{}^{s}
(A^{\point},{B}^{\point,q})}
\quad\hbox{and}\quad
\widetilde\sigma^{s,q}_{\Hom}=\sigma_{\Hom_{}^{s}
(A^{\point},\widetilde{B}^{\point,q})}\; . 
$$ 
}
\proof{
This is a consequence of the previous lemma.}
\smallskip

\prop{NN} 
{
Using the previous notation, 
let  $B_{\tot}^{\point}$ indicate the total complex associated
to ${B}^{\point,\point}$.
Then the complex
$\Hom_{}^{\point}
(A^{\point},B_{\tot}^{\point})$
is canonically identified with: 
\item{$(1)$} 
the total complex of
$\l(\Hom_{}^{s}
(A^{\point},{B}^{p,\point})\r)_{s,p}$: 
$$
\Hom_{}^{\point}
(A^{\point},B_{\tot}^{\point})
\isom
\l(\l(\Hom_{}^{s}
(A^{\point},{B}^{p,\point})\r)_{s,p}\r)_{\tot}
\; ;
$$
\item{$(1')$} 
the total complex of
$\l(\widetilde\Hom_{}^{s}
(A^{\point},\widetilde{B}^{\point,q})\r)_{s,q}$: 
$$
\Hom_{}^{\point}
(A^{\point},B_{\tot}^{\point}) 
\isom
\l(\l(\widetilde\Hom_{}^{s}
(A^{\point},\widetilde{B}^{\point,q})\r)_{s,q}\r)_{\tot}
\; .
$$
}
\proof{ 
The corresponding terms being clearly isomorphic, 
we only have to prove that the differentials in the two
complexes coincide. 
\item{$(1)$} 
Consider the family 
$\Phi=\{\phi_{r}^{a,p}:A^{a}\ra B^{p,a{+}r{-}p}\}_{a,p}$ 
which describes an element of 
$\Hom_{}^{r}(A^{\point},B_{\tot}^{\point})$ 
as well as of 
$\bigoplus_{s{+}p=r}\Hom_{}^{s}(A^{\point},{B}^{p,\point})$.  
Its image using the differential of the first complex is 
$$ 
\eqalign{
D(\Phi)_{r{+}1}^{a,p}
&= 
(d_{B_{\tot}}\Phi)_{r{+}1}^{a,p}+(-1)^{r{+}1} (\Phi d_{A})_{r{+}1}^{a,p}\cr 
&=  
d_{B}^{\prime\,p{-}1,a{+}r{-}p{+}1}\phi_{r}^{a,p{-}1} + 
(-1)^{p}d_{B}^{\prime\prime\,p,a{+}r{-}p}\phi_{r}^{a,p} + 
(-1)^{r{+}1}\phi_{r}^{a{+}1,p}d_{A}^{a}\; . 
}$$ 
Using the second complex and the following diagram 
$$ 
\diagram{
\Hom_{}^{s}(A^{\point},B^{p,\point}) 
&\R\W^{D^{\prime\prime\,s,p}} 
&\Hom_{}^{s{+}1}(A^{\point},{B}^{p,\point}) \cr 
&&\U\H\rn{D^{\prime\,s{+}1,p{-}1}} \cr 
&&\Hom_{}^{s{+}1}(A^{\point},{B}^{p{-}1,\point}),
} 
$$ 
one sees that for $s=r-p$ the image is 
$$ 
\eqalign{
D(\Phi)_{r{+}1}^{a,p}
&= 
D'(\Phi)_{r{+}1}^{a,p} + (-1)^{p} D''(\Phi)_{r{+}1}^{a,p} \cr 
&= 
d_{B}^{\prime\,p{-}1,a{+}r{-}p{+}1}\phi_{r}^{a,p{-}1} + 
(-1)^{p}d_{B}^{\prime\prime\,p,a{+}r{-}p}\phi_{r}^{a,p} + 
(-1)^{r{+}1}\phi_{r}^{a{+}1,p}d_{A}^{a} \; .
}$$ 
So the two differentials coincide. 
\item{$(1')$} 
In this case  the differential of the first complex 
is 
$$ 
\eqalign{
D(\Phi)_{r{+}1}^{a,q}
&= 
(d_{B_{\tot}}\Phi)_{r{+}1}^{a,q}+(-)^{r{+}1} (\Phi d_{A})_{r{+}1}^{a,q}\cr 
&=  
d_{B}^{\prime\,a{+}r{-}q,q}\phi_{r}^{a,q} + 
(-1)^{a{+}r{-}q{+}1}
d_{B}^{\prime\prime\,a{+}r{-}q{+}1,q{-}1}\phi_{r}^{a,q{-}1} + 
(-1)^{r{+}1}\phi_{r}^{a{+}1,q}d_{A}^{a} \; .
}$$ 
Using the second complex the image is 
$$ 
\eqalign{
D(\Phi)_{r{+}1}^{a,q}
&= 
\widetilde D'(\Phi)_{r{+}1}^{a,q} + 
(-1)^{r{-}q{+}1} \widetilde D''(\Phi)_{r{+}1}^{a,q} \cr 
&= 
d_{B}^{\prime\,a{+}r{-}q,q}\phi_{r}^{a,q} + 
(-1)^{a{+}r{-}q{+}1}
d_{B}^{\prime\prime\,a{+}r{-}q{+}1,q{-}1}\phi_{r}^{a,q{-}1} + 
(-1)^{r{+}1}\phi_{r}^{a{+}1,q}d_{A}^{a} \; .
}$$ 
So the two differentials so indeed coincide. 
}

%%%%%%%%%%%%%%%%%%%%%%%%%%%%%%%%%%%%%%%%%%%

\theo{NN} 
{\rif{lmb}
The following sets are canonically isomorphic: 
\item{$(1)$} 
the set 
$\Hom_{} (A^{\point},B_{\tot}^{\point})$
of morphisms of complexes 
between a complex and the total complex of a bicomplex; 
\item{$(2)$} 
the set of cycles 
$Z^{0}\l(\Hom_{}^{\point}(A^{\point},B_{\tot}^{\point})\r)$ ; 
\item{$(3)$} 
the set of cycles 
$Z^{0}\l(\l(\l(\Hom_{}^{s}
(A^{\point},{B}^{p,\point})\r)_{s,p}\r)_{\tot}\r)$ ;
\item{$(3')$} 
the set of cycles 
$Z^{0}\l(\l(\l(\widetilde\Hom_{}^{s}
(A^{\point},\widetilde{B}^{\point,q})\r)_{s,q}\r)_{\tot}\r)$ ;
\item{$(4)$} 
the set of families of maps 
$\{\phi^{\point,p}\in\Hom^{-p}(A^{\point},{B}^{p,\point})\}_{p}$ 
such that 
$$ 
(-1)^{p}d_{\Hom}^{-p}(\phi^{\point,p}) 
+d_{B}^{\prime p{-}1,\point{-}p{+}1} \phi^{\point,p{-}1}
=0
$$ 
for any $p$; 
\item{$(4')$}
the set of families of maps 
$\{\phi^{\point,q}\in
\widetilde\Hom^{-q}(A^{\point},\widetilde{B}^{\point,q})\}_{q}$ 
such that 
$$ 
d_{\widetilde\Hom}^{-q}(\phi^{\point,q}) 
+(-1)^{q{-}1}d_{\widetilde B}^{\prime\prime \point{-}q{+}1,q{-}1} 
\phi^{\point,q{-}1}
=0
$$ 
for any $q$; 
\item{$(5)$} 
the set of families of maps (not maps of complexes in general)
$\{\phi^{\point,p}:A^{\point}\ra {B}^{p,\point}[-p]\}_{p}$ 
such that the ``defect of commutativity'' 
with differentials of one map 
is ``corrected'' by the previous map: 
$$ 
\phi^{\point{+}1,p} d_{A}^{\point}- 
d_{{B}^{p,\point}[-p]}^{\point}\phi^{\point,p}
=
d_{B}^{\prime p{-}1,\point{-}p{+}1} \phi^{\point,p{-}1} 
$$ 
for any $p$; 
\item{$(5')$} 
the set of families of maps (not maps of complexes in general)
$\{\phi^{\point,q}:A^{\point}\ra {B}^{\point,q}[-q]\}_{q}$ 
such that the ``defect of commutativity'' with differentials 
of one map is ``corrected'' by the previous map:
$$ 
\phi^{\point{+}1,q} d_{A}^{\point}- 
d_{\widetilde{B}^{\point,q}[-q]}^{\point}\phi^{\point,q}
=
d_{\widetilde B}^{\prime\prime \point{-}q{+}1,q{-}1} 
\phi^{\point,q{-}1}
$$ 
for any $q$; 
\item{$(6)$} 
 the set of families of maps 
 $\{\phi^{a,p}:A^{a}\ra B^{p,a{-}p}\}_{a,p}$ 
 satisfying the following conditions 
 $$ 
\phi^{a{+}1,p}d_{A}^{a}
-(-1)^{p}d_{B}^{\prime\prime p,a{-}p}\phi^{a,p} 
= d_{B}^{\prime p{-}1,a{-}p{+}1}\phi^{a,p{-}1} 
\leqno{\bf (\NNN)}\rif{conmortot}$$ 
for any $a,p$; 
\item{$(6')$} 
the set of families of maps 
$\{\phi^{a,q}:A^{a}\ra B^{a{-}q,q}\}_{a,q}$ 
satisfying the following conditions 
$$ 
\phi^{a{+}1,q}d_{A}^{a}
-d_{B}^{\prime a{-}q,q}\phi^{a,q} 
= (-1)^{a{-}q{+}1}d_{B}^{\prime\prime a{-}q{+}1,q{-}1}\phi^{a,q{-}1} 
$$ 
for any $a,q$; 
}
\proof{ 
The equivalence of $(1)$ and $(2)$ are well known, 
and the equivalence with $(3)$ and $(3')$ follows 
from (1) and (1${}'$), respectively, of) the Proposition 2.4. 
\endgraf 
The descriptions (6) 
 and (6${}'$) 
follow directly from 
(1) by making  the commutativity condition for 
 the morphism with the differentials explicit. 
\endgraf 
The descriptions $(4)$ and $(4')$ follow directly from 
$(3)$ and $(3')$, respectively. 
\endgraf 
The descriptions $(5)$ and $(5')$ are reformulations of  
$(4)$ and $(4')$, respectively, as well as of 
$(6)$ and $(6')$, respectively. 
}
\smallskip 

%%%%%%%%%%%%%%%%%%%%%%%%%%%%%%%%%%%%%%%%%%%

\rema{NN}
{
The descriptions $(6)$ 
and $(6')$ are useful for applications, 
while $(4)$ and $(4')$ give the most intuitive construction 
of the morphisms from a complex to the total complex associated to a 
bicomplex: we have to define for any $i$ 
a map of graded objects from the complex 
to the $i$-th row (resp. column) shifted by  $-i$,  
in such a way that 
each map is the defect of commutativity for the differentials 
of the next map. 
}
\smallskip 

%%%%%%%%%%%%%%%%%%%%%%%%%%%%%%%%%%%%%%%%%%%%%
\section{Comparison between Saito and HL-localizations.}
%%%%%%%%%%%%%%%%%%%%%%%%%%%%%%%%%%%%%%%%%%%%%

\rema{NN}
{
Let $\c F^\point$ be an object in $C_1(\c O_X,\Diff_X)$.
By definition the differential
$d^i_{\c F}:\c F^i\lra \c F^{i+1}$ is a differential
operator of order one.
So it defines in a unique way a morphism
$\ov d^i_{\c F}: \c F^i\lra \c F^{i+1}\otimes_{\c O_X}\c D_{X,1}$.
%as in \cite{SDDO}.
Locally for each section $s$ of $\c F^i$, 
$\ov d^i_{\c F}(s)$ is a section of
$\c F^{i+1}\otimes_{\c O_X}\c D_{X,1}$
so it may be locally written in a unique way as
$$
\ov d^i_{\c F}(s)=
d^i_{\c F}(s)\otimes 1+\sum_{j=1}^n d^i_{x_j}(s)\otimes
\ds{\partial\over\partial x_j}
\leqno{\bf (\NNN)}\rif{decdo}
$$
using the $\c O_X$-base of $\c D_{X,1}$ given in local coordinates by 
$1, \ds{\partial\over\partial x_1},\cdots , \ds{\partial\over\partial
x_n}$.
The maps $d^i_{\c F}:\c F^i\lra\c F^{i+1}$ are the differentials of the
complex $\c F^\point$;
%(this is a consequence of the commutativity of
%the diagram in \cite{SDDO}).
while $d^i_{x_j}:\c F^i\lra \c F^{i+1}$ are maps
of abelian sheaves.
}
\smallskip

\defi{NN}
{
Let by definition
$\sigma ^{i,j}_{\c F}:\c F^i\lra \c F^{i+j}\otimes_{\c O_X}\cTheta^{-j}_{X}$ 
be the maps defined as follow: 
$$
\xymatrix{
\c F^{i}  \ar[r]^(.3){\ov{d}}\ar[drrr]_{\sigma^{i,j}_{\c F}} 	& 
\c F^{i+1}\otimes_{\c O_X}\c D_{X,1}
\ar[r]^(.7){\ov{d}\otimes{id}}				& 
\cdots\ar[r]						&
\c F^{i+j}\otimes_{\c O_X}\c D_{X,1}\otimes_{\c O_X}\cdots
\otimes_{\c O_X}\c D_{X,1}\ar[d] 			\\
							&
							&
							&
\c F^{i+j}\otimes_{\c O_X}\cTheta^{-j}_{X}.		\\ 
}
$$
The vertical map is the identity on $\c F^{i+1}$ tensor the map
obtained by the composition of the projections  
$\c D_{X,1}\lra \cTheta^1_X$ and 
$\cTheta^1_X \otimes_{\c O_X}\cdots\otimes_{\c O_X}\cTheta^1_X
\lra \cTheta^{-j}_{X}$.
We observe that these maps $\sigma^{i,j}_{\c F}$ are related to the
structural  morphisms of the $\cOmega^\point_{X}$-module 
$$
%m^{i,j}_{\c F}:
\c F^i\otimes_{\c O_X}\cOmega^j_{X}\lra \c F^{i+j};$$
in fact they are adjoints
%(by (\cite{FDO}) and 
because $\cTheta^{-j}_X\bydef\wedge^j\cTheta_{X}\cong
\cHom_{\c O_X}(\cOmega^j_{X},\c O_X)$)
and
$$
\eqalign{
\Hom_{\c O_X}(\c F^i\otimes_{\c O_X}\cOmega^i_X,\c F^{i+j})
&\cong  \Hom_{\c O_X}(\c F^i,\cHom_{\c O_X}(\cOmega^i_X,\c F^{i+j}))\cr
&\cong  \Hom_{\c O_X}(\c F^i,\c F^{i+j}\otimes_{\c O_X}\cTheta^{-j}_{X}).\cr
}
$$
}
\smallskip

\defi{NN}
{
Let $\c F^{\point}\in C_1(\c O_X,\Diff_X)$. We define for each
$i\in \Bbb Z$ and $j\in \{0,...,n\}$ ($n=\dim X$) 
the maps
$$
\eta^{i,0}_{\c F}=id_{\c F^i}: \c F^i \lra \c F^i
\leqno{\bf (\NNN)}\rif{eta0}
$$
and
$$\matrix{
\eta^{i,j}_{\c F}: & \c F^i & \lra &
\c F^{i+j}\otimes_{\c O_X}\cTheta^{-j}_{X}\hfill\cr
&
s &\longmapsto &
\sum\limits_{i_1<... <i_j}d^{i+j-1}_{x_{i_j}}\opp\cdots\opp 
d^{i}_{x_{i_1}}(s)\otimes\ds{\partial\over{\partial x_{i_1}}}\wedge
\cdots\wedge {\partial\over {\partial x_{i_j}}}\cr
} 
\leqno{\bf (\NNN)}\rif{etaj}
$$
for $j\in \{1,...,n\}$.
Then (up to a sign $(-1)^{j\choose 2}$)  we have that  
$\sigma^{i,j}_{\c F}=j!\eta^{i,j}_{\c F}$,
so these maps $\eta^{i,j}_{\c F}$ do not depend on local coordinates. 
}
\smallskip

%%%%%%%%%%%%%%%%%%%%%%%%%%%%%%%%%%%%%%%%%%%%%%%
%%%%%%%%%% Lemma condizione d d_{x_i} %%%%%%%%%%%%%%%%%%%%%%
%%%%%%%%%%%%%%%%%%%%%%%%%%%%%%%%%%%%%%%%%%%%%%%

\lemm{NN}
{\rif{condiff}
Given $\c F^\point\in C_1(\c O_X,\Diff_X)$; the morphisms $d^i_{\c F}$, 
$d^i_{x_j}$ of (\cite{decdo}) for $i\in \Bbb Z$ and $j\in \{0,...,n\}$
satisfy the following conditions:
	\item{i)} $d^{i+1}_{\c F}\opp d^{i}_{\c F}=0$
	\item{ii)} $d^{i+1}_{x_j}\opp d^i_{\c F}+
	            d^{i+1}_{\c F}\opp d^i_{x_j}=0$
	\item{iii)}$d^{i+1}_{x_j}\opp d^i_{x_k}+
		    d^{i+1}_{x_k}\opp d^i_{x_j}=0$	
	\item{iv)} $d^{i+1}_{x_j}\opp d^i_{x_j}=0$.
}
\proof{
The first condition is given by the hypothesis $\c F^\point\in
C_1(\c O_X,\Diff_X)$.
The conditions ii) to iv) follow from the condition that the
composition
$$
\ov{d^{i+1}_{\c F}\opp d^i_{\c F}}:
\c F^i
\lra 
\c F^{i+2}\otimes_{\c O_X}\c D_{X,2}
$$
is zero because it  corresponds to 
$d^{i+1}_{\c F}\opp d^i_{\c F}=0$.
}
\smallskip

%%%%%%%%%%%%%%%%%%%%%%%%%%%%%%%%%%%%%%%%%
%%%%%%%%% definizione Phi	%%%%%%%%%%%%%%%%%%%%%%%
%%%%%%%%%%%%%%%%%%%%%%%%%%%%%%%%%%%%%%%%%

\defi{NN}
{\rif{defiphi}
Let $i:\c O_X\lra \c D_{X}$ be the usual inclusion which is linear for
both the $\c O_X$-Module structures of $\c D_{X}$.
Given $\c F^{\point}\in C_1(\c O_X,\Diff_X)$
we define the morphisms
$$
\Phi^i_{\c F}: \c F^i\lra \bigoplus_{j=0}^n 
\c F^{i+j}\otimes_{\c O_X}\c D_{X}\otimes_{\c O_X}\cTheta^{-j}_{X}
$$
for each $i\in \Bbb Z$
in the following way:
we consider the composition
$$
\xymatrix{
\c F^i\ar[r]^(.4){\eta^{i,j}_{\c F}}	\ar[dr]_(.4){\Phi^{i,j}_{\c F}}	&
\c F^{i+j}\otimes_{\c O_X}\cTheta^{-j}_{X}
\ar[d]^{id_{\sss{\c F^{i+j}}}\otimes i\otimes id_{\sss{\cTheta^{-j}}}}	\\
									&
\c F^{i+j}\otimes_{\c O_X}\c D_{X}\otimes_{\c O_X}\cTheta^{-j}_{X}	\\
}
$$
and by definition
$\Phi^i_{\c F}:=\sum_{j=0}^n \Phi^{i,j}_{\c F}$.
}
\smallskip

We want to prove that the morphisms 
$\Phi^i_{\c F}:\c F^i\lra ({\tDR}_{X}{\tDR}^{-1}_{X}
(\c F^\point))^i$  
define a morphism of complexes. 
%In order to do this we need a general result about complexes.
\smallskip

%%%%%%%%%%%%%%%%%%%%%%%%%%%%%%%%%%%%%%%%
%%%%%% Teorema commutativita' Phi differenziali		         %%%%%%%
%%%%%%%%%%%%%%%%%%%%%%%%%%%%%%%%%%%%%%%%

\theo{NN}
{
The maps $\Phi^i_{\c F}$ of \cite{defiphi}
define 
$$
\Phi_{\c F}:
\c F^{\point}\lra {\tDR_X}{\tDR_X}^{-1}(\c F^{\point})
$$
which is a morphism of complexes in $C_1(\c O_X,\Diff_X)$.
}
%%%%%%%%%%%%%%%%%%%%%%%%%%%%%%%%%%%%%%%%
%%%% Dimostrazione teorema commutativita' Phi d %%%%%%%%%%%%%%
%%%%%%%%%%%%%%%%%%%%%%%%%%%%%%%%%%%%%%%%
\proof{
We have to prove that
the diagram
$$
\xymatrix{
\c F^i	
\ar[r]^{\sss{d_{\c F}^i}}	\ar[d]^{\Phi^i_{\c F}}	&
\c F^{i+1}\ar[d]^{\Phi^{i+1}_{\c F}}		\\
\bigoplus_{j=0}^n 
\c F^{i+j}\otimes_{\c O_X}\c D_{X}
\otimes_{\c O_X}\cTheta^{-j}_{X}
\ar[r]_{d_{\sss{{\tDR}{\tDR}^{-1}}}^i}&
\bigoplus_{j=0}^n 
\c F^{i+1+j}\otimes_{\c O_X}\c D_{X}
\otimes_{\c O_X}\cTheta^{-j}_{X}\\
}
$$
is commutative.
\endgraf
We recall that 
$
{\tDR_X}{\tDR_X}^{-1}(\c F^{\point})=
(\c G^{\point\point})_{tot}
$
where
$$
\c G^{p,q}=
\c F^q\otimes_{\c O_X}\c D_{X}\otimes_{\c O_X}\cTheta^p_{X}
$$
and
$$
d^{\prime p,q}_{\c G}:
\c F^q\otimes_{\c O_X}\c D_{X}\otimes_{\c O_X}\cTheta^p_{X}
\lra  
\c F^{q}\otimes_{\c O_X}\c D_{X}\otimes_{\c O_X}\cTheta^{p+1}_{X}
$$
is
$$
d^{\prime p,q}_{\c G}=	
id_{\c F^q}\otimes d^{p}_{Sp\c O_X} ; 	
\leqno{\bf (\NNN)}\rif{d2}
$$
while
$$
d^{\prime\prime p,q}_{\c G}:
\c F^q\otimes_{\c O_X}\c D_{X}\otimes_{\c O_X}\cTheta^p_{X}
\lra  
\c F^{q+1}\otimes_{\c O_X}\c D_{X}\otimes_{\c O_X}\cTheta^p_{X}
$$ 
is 
$$
d^{\prime\prime p,q}_{\c G}=					
{\tDR_X}^{-1}(d^q_{\c F})\otimes 
id_{\cTheta^p_{X}}\; .
\leqno{\bf (\NNN)}\rif{d1}
$$
%where 
%$d^{j}_{Sp\c O_X}$ 
%was defined in (\cite{diffSPOX}).
\endgraf

%%%%%%%%%%%%%%%%%%%%%%%%%%%%%%%%%%%%%%%%%%%%%%%%%%%%%
%%%% citazione lemma 					%%%%%%%%%%%%%%%%%%%%%%%%%%%%%%
%%%%%%%%%%%%%%%%%%%%%%%%%%%%%%%%%%%%%%%%%%%%%%%%%%%%%
Now by Lemma \cite{lmb}, we have only to prove that
(\cite{conmortot}) 
$$
\Phi^{a{+}1,p}_{\c F}\opp d^a_{\c F}-
(-1)^p d^{\prime\prime {-}p,a{+}p}_{\c G}
\opp \Phi^{a,p}_{\c F}=
d^{\prime {-}p{-}1,a{+}p{+}1}_{\c G}\opp \Phi^{a,p{+}1}_{\c F}
$$
is true.\endgraf
Let $s$ be a section of $\c F^a$, then
%%%%%%%%%%%%%%%%%%%%%%%%%%%%%%%%%%%%%%%%%%%%%%%%%%%%%
%%%% phi d_{\c F}	 %%%%%%%%%%%%%%%%%%%%%%%%%%%%
%%%%%%%%%%%%%%%%%%%%%%%%%%%%%%%%%%%%%%%%%%%%%%%%%%%%%
$$
\Phi^{a{+}1,p}_{\c F}\opp d^a_{\c F}(s)=	 		
\sum_{i_1<\cdots <i_p}
d^{a+p}_{x_{i_p}}\cdots d^{a+1}_{x_{i_1}}d^a_{\c F}(s)
\otimes 1 \otimes
\ds{\partial\over\partial x_{i_1}}\wedge \cdots
\wedge \ds{\partial\over\partial x_{i_p}};
%\leqno{\bf (\NNN)}\rif{phid}
$$
while
%%%%%%%%%%%%%%%%%%%%%%%%%%%%%%%%%%%%%%%%%%%%%%%%%%%%%
%%%% d_1 phi 		 %%%%%%%%%%%%%%%%%%%%%%%%%%%%
%%%%%%%%%%%%%%%%%%%%%%%%%%%%%%%%%%%%%%%%%%%%%%%%%%%%%
$$
\eqalign
{
& d^{\prime\prime {-}p,a{+}p}_{\c G}\opp \Phi^{a,p}_{\c F}(s)=	\cr
&= d^{\prime\prime {-}p,a{+}p}_{\c G}
\big(
\sum_{i_1<\cdots <i_p}
d^{a+p-1}_{x_{i_p}}\cdots d^{a}_{x_{i_1}}(s)
\otimes 1 \otimes
\ds{\partial\over\partial x_{i_1}}\wedge \cdots
\wedge \ds{\partial\over\partial x_{i_p}}\big)=	\cr
%%%%%%%%%%%%%%%%%%%%%%%%%%%%%%%%%%%%%%%%%%%%%%%%%%%%%
&  =						
\sum_{i_1<\cdots <i_p}
d^{a+p}_{\c F}\big(
d^{a+p-1}_{x_{i_p}}\cdots d^{a}_{x_{i_1}}(s)
\big)
\otimes 1 \otimes
\ds{\partial\over\partial x_{i_1}}\wedge \cdots
\wedge \ds{\partial\over\partial x_{i_p}}+
						\cr
%%%%%%%%%%%%%%%%%%%%%%%%%%%%%%%%%%%%%%%%%%%%%%%%%%%%%
&
+
\sum_k\sum_{i_1<\cdots <i_p}
d^{a+p}_{x_k}
\big(
d^{a+p-1}_{x_{i_p}}\cdots d^{a}_{x_{i_1}}(s)
\big) 
\otimes \ds{\partial\over\partial x_{k}}\otimes 
\big(
\ds{\partial\over\partial x_{i_1}}
\wedge\cdots \wedge
\ds{\partial\over\partial x_{i_p}}
\big)=						\cr
%%%%%%%%%%%%%%%%%%%%%%%%%%%%%%%%%%%%%%%%%%%%%%%%%%%%%
&    =
\sum_{i_1<\cdots <i_p}(-1)^p
d^{a+p}_{x_{i_p}}\cdots d^{a+1}_{x_{i_1}}d^a_{\c F}(s)
\otimes 1 \otimes
\ds{\partial\over\partial x_{i_1}}\wedge \cdots
\wedge \ds{\partial\over\partial x_{i_p}}+
						\cr
%%%%%%%%%%%%%%%%%%%%%%%%%%%%%%%%%%%%%%%%%%%%%%%%%%%%%
&
+
\sum_{k}
\sum_{i_1<\cdots <i_{p+1}}(-1)^{p+k}
d^{a+p}_{x_{i_{p+1}}}\cdots d^{a}_{x_{i_{1}}}(s)
\otimes 
\ds{\partial\over\partial x_{i_k}}\otimes 
\big(
\ds{\partial\over\partial x_{i_1}}
\wedge\cdots \wedge
\widehat{\ds{\partial\over\partial x_{i_k}}}
\wedge\cdots \wedge
\ds{\partial\over\partial x_{i_{p+1}}}
\big).\cr}
%\leqno{\bf (\NNN)}\rif{d1phi}
$$
On the right hand  side of 
(\cite{conmortot}) 
we obtain
%%%%%%%%%%%%%%%%%%%%%%%%%%%%%%%%%%%%%%%%%%%%%%%%%%%%%%%%%%%%
%%%%% d_2 Phi %%%%%%%%%%%%%%%%%%%%%%%%%%%%%%%%%%%%%%%%%%%%%%
%%%%%%%%%%%%%%%%%%%%%%%%%%%%%%%%%%%%%%%%%%%%%%%%%%%%%%%%%%%%
$$
\eqalign
{&
d^{\prime {-}p{-}1,a{+}p{+}1}_{\c G}\opp \Phi^{a,p{+}1}_{\c F}(s)=\cr
&=d^{\prime {-}p{-}1,a{+}p{+}1}_{\c G}
\big(
\sum_{i_1<\cdots <i_{p+1}}
d^{a+p}_{x_{i_{p+1}}}\cdots d^a_{x_{i_{1}}}(s)
\otimes 1 \otimes
\ds{\partial\over\partial x_{i_1}}
\wedge\cdots  \wedge
\ds{\partial\over\partial x_{i_{p+1}}}
\big)=						\cr
&   =
\sum_{i_1<\cdots <i_{p+1}}
\sum_{k}(-1)^{k+1}
d^{a+p}_{x_{i_{p+1}}}\cdots d^a_{x_{i_{1}}}(s)
\otimes 
\ds{\partial\over\partial x_{i_k}}\otimes 
\big(
\ds{\partial\over\partial x_{i_1}}
\wedge\cdots \wedge
\widehat{\ds{\partial\over\partial x_{i_k}}}
\wedge\cdots \wedge
\ds{\partial\over\partial x_{i_{p+1}}}
\big).						\cr
}
$$
Thus
we have established our assertion.
}
\smallskip

%%%%%%%%%%%%%%%%%%%%%%%%%%%%%%%%%%%%%%%%%%%%%%%%%%%%%%%%%%%%
%%%% Teorema di confronto Saito HL     %%%%%%%%%%%%%%%%%%%%%
%%%%%%%%%%%%%%%%%%%%%%%%%%%%%%%%%%%%%%%%%%%%%%%%%%%%%%%%%%%%

\theo{NN}
{
The functor
$$
\lambda_1: D_1(\c O_X,\Diff_X)\lra D(\c O_X,\Diff_X)
$$
is an equivalence of categories with quasi-inverse the functor
$$
{\tDR}_{1,X}{\tDR}^{-1}_{X}:
D(\c O_X,\Diff_X)\lra D_1(\c O_X,\Diff_X).
$$
}
%%%%Dimostrazione
\proof{
Let $G:={\tDR}_{1,X}\circ {\tDR}^{-1}_{X}$.
\endgraf
By Saito's results we obtain
$\lambda_1\opp G={\tDR}_{X}\circ
{\tDR}^{-1}_{X}
\buildrel{\cong}\over\lra
id_{D(\c O_X,\Diff_X)}$.
%\endgraf
We want to prove that there exists an isomorphism of functors
$id_{D_1(\c O_X,\Diff_X)}\lra G\opp \lambda_1$.
\endgraf
In \cite{defiphi} we defined a functorial morphism 
$\Phi^{\point}_{\c F}:\c F^\point\lra G\opp \lambda_1(\c F^\point)$
(for each $\c F^\point\in D_1(\c O_X,\Diff_X)$) which 
is a ${\tDR}^{-1}_{X}$-quasi-isomorphism
because the triangle
$$
\xymatrix{
\c F^\point\ar[rd]_{id_{\c F^\point}}
\ar[r]^(.4){\Phi^{\point}_{\c F}}						&
G\opp \lambda_1(\c F^\point)	
\ar[d]								\\
								&
\c F^\point							\\
}$$
commutes (where the vertical map is that
induced by the projection $\c D_{X}\lra \c O_X$
and it is 
a ${\tDR}^{-1}_{X}$-quasi-isomorphism by Saito's result).
So the morphism defined by the functor $\Phi:id_{D_1(\c O_X,\Diff_X)}\lra G\opp \lambda_1$
is an isomorphism.
}
\eject
%%%%%%%%%%%%%%%%%%%%%%%%%%%%%%%%%%%%%%%%%%%%%%%%%%%%%%%%%%%%%%%%%%%%%%%%%%%%%
\section{Comparison between Saito and Du Bois categories.}
%%%%%%%%%%%%%%%%%%%%%%%%%%%%%%%%%%%%%%%%%%%%%%%%%%%%%%%%%%%%%%%%%%%%%%%%%%%%%

\defi{NN}
{\rif{CDIFFDB}
By definition the category $CF_1(\c O_X,\Diff_X)$ is the category whose
objects are filtered complexes 
$(K^{\point}, d, F)$ 
such that:
 \item{i)} $K^{\point}$ is a complex of $\c O_X$-Modules;
 \item{ii)} $F$ is a decreasing filtration on $K^{\point}$ given by
	sub-$\c O_X$-Modules and 
	$F$ is biregular 
	(that is on every component $K^i$ of $K^{\point}$, 
	$F$ induces a finite filtration; 
	so there exist integers $p$ and $q$ such that $F^p K^i=K^i$ and $F^q K^i=0$);
 \item{iii)} $d$ is a relative differential operator of order at most
	1 which respects the filtrations;
 \item{iv)} $\gr_F(d)$ is $\c O_X$-linear.
\endgraf
Morphisms in $CF_1(\c O_X,\Diff_X)$ are 
$f^{\point}: \c F^{\point}\lra\c G^{\point}$
$\c O_X$-linear maps commuting with differentials and compatible with
fitrations [DB.1, 1].
}
\smallskip

%%%%%%%%%%%%%%%%%%%%%%%%%%%%%%%%%%%%%%%%%%%%%%%%%%%%%%%%%%%%%%%%%%%%%%%%%%%%%

\rema{NN}
{
The complex $\cOmega^{\point}_{X}$ is filtered 
by truncation so 
$Gr^p_F(\cOmega^{\point}_{X})=\cOmega^p_{X}$.
} 
\smallskip

%%%%%%%%%%%%%%%%%%%%%%%%%%%%%%%%%%%%%%%%%%%%%%%%%%%%%%%%%%%%%%%%%%%%%%%%%%%%%

\defi{NN}
{ 
A filtered homotopy between 
$u,v: \c F^{\point}\lra\c G^{\point}$ 
is a homotopy $h$ such that 
$h^i(F^p\c F^i)\subseteq F^p\c G^{i-1}$ 
and $h^i$ is $\c O_X$-linear.
\endgraf
Let $KF_1(\c O_X,\Diff_X)$ be the category whose objects are those of
$CF_1(\c O_X,\Diff_X)$ and whose morphisms are the classes of morphisms in
$CF_1(\c O_X,\Diff_X)$ up to homotopy.
}
\smallskip

%%%%%%%%%%%%%%%%%%%%%%%%%%%%%%%%%%%%%%%%%%%%%%%%%%%%%%%%%%%%%%%%%%%%%%%%%%%%%

\defi{NN}
{\rif{DDIFFDB}
A morphism $f^{\point}: \c F^{\point}\lra\c G^{\point}$
is said a filtered quasi-isomorphism if $\gr f$ is a quasi-isomorphism
where $\gr$ is the functor
$$
\gr : CF_1(\c O_X,\Diff_X)\lra CG(X)
$$
sending a filtered complex $(\c F^\point,F)$ to its graded complex.
(Here $CG(X)$ is the category of complexes of $\c O_X$-Modules with a
finite graduation in each degree.)
\endgraf
Let $DF_1(\c O_X,\Diff_X)$ be the category obtained by localizing $KF_1(\c O_X,\Diff_X)$ with
respect to filtered quasi-isomorphisms.
}
\smallskip

%%%%%%%%%%%%%%%%%%%%%%%%%%%%%%%%%%%%%%%%%%%%%%%%%%%%%%%%%%%%%%%%%%%%%%%%%%%%%

\rema {NN}
{ 
The categories $KF_1(\c O_X,\Diff_X)$ and $DF_1(\c O_X,\Diff_X)$ are triangulated
categories where 
\endgraf
\item{i)} the shift functor is the usual one: 
$T(K)^i=K^{i+1}$,
$d_{T(K)}=-d_K$ and $F_{T(K)}=F_K$;
\item{ii)} if $f^{\point}: \c F^{\point}\lra\c G^{\point}$ is a morphism in
$CF_1(\c O_X,\Diff_X)$ its mapping cone is 
$\c M^{\point}\bydef T(\c
F^{\point})\oplus \c G^{\point}\in C_1(\c O_X,\Diff_X)$ 
with filtration defined by
$F^p(\c M^{\point})= F^p T(\c F^{\point})\oplus F^n
\c G^{\point}$
and
differential $d_{\c M}$ defined by the matrix
$\l[
\matrix{
d_{T(\c F)} & 0 \cr
T(f)& d_{\c G}\cr
}\r]$.
}
\smallskip
%%%%%%%%%%%%%%%%%%%%%%%%%%%%%%%%%%%%%%%%%%%%%%%%%%%%%%%%%%%%%%%%%%%%%%%%%%%%%%%%
%%%%%%%%%%%%%%%%%%%%%%%%%%%%%%%%%%%%%%%%%%%%%%%%%%%%%%%%%%%%%%%%%%%%%%%%%%%%%%%%
%%%%%%%%%%%%%%%%%%%%%%%%%%%%%%%%%%%%%%%%%%%%%%%%%%%%%%%%%%%%%%%%%%%%%%%%%%%%%%%%
%%%%%%%%%%%%%%%%%%%%%%%%%%%%%%%%%%%%%%%%%%%%%%%%%%%%%%%%%%%%%%%%%%%%%%%%%%%%%%%%

\defi{NN}
{
Let $(\c L_j,F)$ with $j\in \{1,2\}$ be two filtered $\c O_X$-modules  
with increasing filtration such that $F_p\c L_j=0$ for $p\ll 0$.
By definition
$$\cHom_{\Diff_{X/S}}((\c L_1,F),(\c L_2,F))$$
is the sheaf of filtered differential operators, that is
$\Phi\in \cHom_{\Diff_{X/S}}((\c L_1,F),(\c L_2,F))$ if and only if
the composition
$$
F_{-p}\c L_1\lra \c L_1
\buildrel{\Phi}\over{\lra} \c L_2\lra
\c L_2/F_{p-q-1}\c L_2\leqno{\bf (\NNN)}\rif{fildiffcon}
$$
 has order at most $q$  for each $p$, $q$.
\endgraf
This condition implies that $\Phi(F_p\c L_1)\subset F_p\c L_2$ and
that the map between the graded objects
$Gr_p^F(\c L_1)\lra Gr_p^F(\c L_2)$ 
is $\c O_X$-linear.
} 
\smallskip

%%%%%%%%%%%%%%%%%%%%%%%%%%%%%%%%%%%%%%%%%%%%%%%%%%%%%%%%%%%%%%%%%%%%%%%%%%%%%

\rema {NN}
{
We observe that if $(\c L_1,F)\lra (\c L_2,F)$
is a filtered differential operator of order one then the condition
(\cite{fildiffcon}) is equivalent to the condition given by Du Bois of
having graded $\c O_X$-linear for the objects $(\c L_j,F^\prime)$ where
$F^{\prime p}\c L_j=F_{-p}\c L_j$ is the opposite filtration.
}
\smallskip

%%%%%%%%%%%%%%%%%%%%%%%%%%%%%%%%%%%%%%%%%%%%%%%%%%%%%%%%%%%%%%%%%%%%%%%%%%%%%

\defi{NN}
{
We denote by $CF(\c O_X,\Diff_X)$ the category of complexes of 
(increasing) filtered $\c O_X$-Modules
and filtered differential operators.
\endgraf
The functor $\lambda_1$ induces a functor which we denote by
$$\lambda_1 F:CF_1(\c O_X,\Diff_X)\lra CF(\c O_X,\Diff_X).$$
It sends objects of $CF_1(\c O_X,\Diff_X)$ into themselves with 
the opposite filtration (which becomes increasing). 
We denote by $DF(\c O_X,\Diff_X)$ the category obtained
localizing $CF(\c O_X,\Diff_X)$ with respect to filtered quasi-isomorphisms.
}
\smallskip

%%%%%%%%%%%%%%%%%%%%%%%%%%%%%%%%%%%%%%%%%%%%%%%%%%%%%%%%%%%%%%%%%%%%%%%%%%%%%

\theo{NN}
{
The functor
$$
\lambda_1F:
DF_1(\c O_X,\Diff_X)\lra DF(\c O_X,\Diff_X)
$$
is an equivalence of categories;
so
also is the functor
$$
{\tDR}_{1,X}F:
DF(\c D_{X})^\sr\lra
DF_1(\c O_X,\Diff_X).
$$
}
\proof{
Saito proved that the functor 
${\tDR}_{X}F:
DF(\c D_{X})^\sr\lra DF(\c O_X,\Diff_X)$ is an equivalence of categories
with quasi-inverse the functor ${\tDR}^{-1}_{X}F$.
We want to extend the result of the previous section to the filtered context.
As for the non filtered case,
let $GF= {\tDR}_{1,X}F\opp{\tDR}^{-1}_{X}F$
we have only to
prove that there is an isomorphism of functors 
$id_{DF_1(\c O_X,\Diff_X)}\lra GF\opp \lambda_1F$.
\endgraf
In the previous section we defined an isomorphism of functors
$$\Phi:
id_{D_1(\c O_X,\Diff_X)}
\lra 
{\tDR}_{1,X}\opp{\tDR}^{-1}_{X}\opp \lambda_1F
$$
into $D_1(\c O_X,\Diff_X)$.
Now given $(\c F^\point,F)\in CF_1(\c O_X,\Diff_X)$ we want to prove that 
the morphism $\Phi^\point_{\c F}$ respects the filtrations, so it induces an
isomorphism of functors also in the filtered case.
We observe that the map $\eta^{i,j}_{\c F}$ satisfies:
$$
\eta^{i,j}_{\c F}:F_p(\c F^i)\lra
F_{p}(\c F^{i+j})\otimes_{\c O_X}\cTheta^{-j}_{X}
\leqno{\bf (\NNN)}\rif{respfil}
$$
because the differentials of the complex respect the filtrations.
\endgraf
We recall that the filtration on the complex 
${\tDR}_{1,X}{\tDR}^{-1}_{X}(\c F^\point)$
is built as explained in [S.1; 2.1.3, 2.1.5] for a single $\c O_X$-Module
$\c F$.
Generalizing this construction to complexes we have that 
$$
F_p(\bigoplus_{j=0}^n \c F^{i+j}\otimes \c D_{X}\otimes \cTheta^j_{X})=
\bigoplus_{j=0}^n \sum_{l\ge 0} F_{p-l}(\c F^{i+j})\otimes\c D_{X,l}\otimes \cTheta^{-j}_{X}.$$
This implies that  also $\Phi^{i,j}_{\c F}$ respects the filtrations because $F_p(\c F^i)$
takes image into
$F_{p}(\c F^{i+j})\otimes\c D_{X,0}\otimes \cTheta^{-j}_{X}$; so we have established our thesis.
}
\smallskip

%%\rema{NN}
%%{
%%We observe that each complex $\c F^\point\in C_1(\c O_X,\Diff_X)$
%%admits a canonical filtration $F$ as follow:
%%$F_c^{-p}\c F^i={\rm Ker}(\eta^{i,p+1}_{\c F}:\c F^i\lra \c F^{i+p+1}\otimes
%%\cTheta^{-(p+1)}_{X})$.
%%In this way we can build a functor:
%%$$
%%\matrix{
%%{\rm Can}:&
%%C_1(\c O_X,\Diff_X)&\lra& CF_1(\c O_X,\Diff_X)\cr
%%&\hfill\c F^\point&\longmapsto&(\c F^\point, F_c)\cr}.
%%\leqno{\bf (\NNN)}\rif{can}
%%$$
%%I do not know if this functor localizes with respect to
%%$\tDR^{-1}_X$-quasi-isomorphism from one side and graded-quasi-isomorphism from the other.
%%}
%%\smallskip

\chiudiriferimenti

\refer

\biblio{B}
                    {Berthelot  P.
                     {Cohomologie cristalline des sch\'emas de
			caract\'eristique $p>0$.}
                     {Lecture Notes in Mathematics, Vol. 407.}
                     {\it Springer-Verlag, Berlin-New York,} 1974.}

\biblio{Bo}
                    {Borel et al.
                     {\it Algebraic D-Modules.}
                     Perspectives in Mathematics, Vol. 2
			J. Coates and S. Helgasan editors.}

\biblio{DB.1}
                    {Du Bois Ph.
                     {Complexe de de Rham filtr\'e d'une vari\'et\'e
			singuli\`ere.}
                     {\it Bull. Soc. Math. France},
                     {\bf 109} (1981), 41--81.}

\biblio{DB.2}
                    {Du Bois Ph.
                     {Dualit\'e dans la cat\'egorie des complexes filtr\'es
			d'op\'erateurs differentiels d'ordre $\leq 1$.}
                     {\it Collect. Math.},
                     {\bf 41} (1990), 89--121.}

\biblio{EGA}
                    {Grothendieck A. and Dieudonn\'e J.
                     El\'ements de g\'eom\'etrie alg\'ebrique.
                     {\it Inst. Hautes \'Etudes Sci. Publ. Math.},
                     {\bf 4} (1960), {\bf 8} (1961), {\bf 11} (1961), {\bf
			17} (1963),
                     {\bf 20} (1964), {\bf 24} (1965), {\bf 28} (1966), {\bf
			32} (1967).}

\biblio{HL}
                    {Herrera M. and Lieberman D.
                     {Duality and the de Rham cohomology of infinitesimal
			neighborhoods.}
                     {\it Invent. Math.},
                     {\bf 13} (1971), 97--124.}

\biblio{M}
                    {Mebkhout Z.
                     {\it Le formalisme des six op\'erations de Grothendieck
			pour les $\c D\sb X$-modules coh\'erents.}
                     {Travaux en Cours, 35},
                     {\it Hermann, Paris,} 1989.}

\biblio{S.1}
                    {Saito M.
                     {Modules de Hodge polarisables.}
                     {\it Publ. RIMS, Kyoto Univ.},
                     {\bf 24} (1988), 849--995.}

\biblio{S.2}
                    {Saito M.
                     {Induced $\c D$-modules and differential complexes.}
                     {\it Bull. Soc. Math. France},
                     {\bf 117} (1989), 361--387.}

\bigskip
\bigskip

\halign{
#\hfill &\qquad#\hfill \cr 
Luisa Fiorot\cr 
Universit\`a degli Studi di Padova &Universit\'e Louis Pasteur, Strasbourg\cr 
Dipartimento di Matematica Pura e Appl.%icata 
&UFR de Math\'ematique \cr 
via Belzoni 7 &7 rue Ren\'e Descartes \cr 
35131 Padova &67084 Strasbourg \cr 
Italy &France \cr 
\tt fiorot@math.unipd.it &\tt fiorot@math.u-strasbg.fr\cr 
}

\end 

\smallskip

Luisa Fiorot

Universit\`a degli Studi di Padova

Dipartimento di Matematica Pura ed Applicata

via Belzoni 7

35131 Padova

Italy

fiorot@math.unipd.it

\smallskip

Universit\'e Louis Pasteur, Strasbourg

UFR de Math\'ematique

7 rue Ren\'e Descartes

67084 Strasbourg

France

fiorot@math.u-strasbg.fr

\end
\lemm{NN}
{The forgetful functor
$$
CF_1(\c O_X,\Diff_X)\lra C_1(\c O_X,\Diff_X)
$$
sends filtered quasi-isomorphisms into
${\tDR}^{-1}_{X}$-quasi-isomorphisms
so it defines a functor
$$
DF_1(\c O_X,\Diff_X)\lra D_1(\c O_X,\Diff_X).
$$
}  
\proof{
Given $\c F^\point\lra \c G^\point$ a morphism in $CF_1(\c O_X,\Diff_X)$ 
which is a filtered quasi-isomorphism, we have that also
${\tDR}^{-1}_{X}(\c F^\point)\lra 
{\tDR}^{-1}_{X}(\c G^\point)$ is a filtered quasi-isomorphism; this
implies that it is a quasi-isomorphism in the category $D(\c D_{X}^\sr)$.
So $\c F^\point$ and $\c G^\point$ are $\tDR^{-1}$-quasi-isomorphism.
}
\smallskip

%%%%%%%%%%%%%%%%%%%%%%%%%%%%%%%%%%%%%%%%%%%%%%%%%%%%%%%%%%%%%%%%%%%%%%
%%%%% Lemma morfismo verso un tot %%%%%%%%%%%%%%%%%%%%%%%%%%%%%%%%%%%%
%%%%%%%%%%%%%%%%%%%%%%%%%%%%%%%%%%%%%%%%%%%%%%%%%%%%%%%%%%%%%%%%%%%%%%
\lemm{NN}
{\rif{lmb}
Let $\c A$ be an abelian category.
For each complex $A^\point$, and bicomplex $B^{\point\point}$ 
bounded
in one degree;
the following are equivalent:
	\item {i)} the data of a morphism of complexes 
	 $\Phi: A^\point\lra (B^{\point\point})_{tot}$ 
	\item {ii)} the data of a family of morphisms
		$\Phi^i_j:A^i\lra B^{i+j,-j}$
		such that
$$
\Phi^{i+1}_j\opp d^i_A
-d^{i+j+1,-(j+1)}_2
\opp \Phi^i_{j+1}=
(-1)^{j}
d^{i+j,-j}_1\opp \Phi^i_{j}
\leqno{\bf (\NNN)}\rif{conmortot}
$$
where
$d^{k,l}_1: B^{k,l}\lra B^{k+1,l}$ is the first differential of the
bicomplex $B^{\point\point}$;
while $d^{k,l}_2: B^{k,l}\lra B^{k,l+1}$ is the second one
($d^{k+1,l}_2\opp d^{k,l}_1=d^{k,l+1}_1\opp d^{k,l}_2$).
\endgraf
Moreover the correspondence is $\Phi^i=\sum_j\Phi^i_j:A^i\lra 
(B^{\point\point})_{tot}^i$.
}
\proof
{
The complex $(B^{\point\point})_{tot}$ is by definition:
$$
(B^{\point\point})_{tot}^i=\bigoplus_{j\in\Bbb Z}B^{i+j,-j}
$$
while
$$
d^i_{tot}: 
\bigoplus_{j\in\Bbb Z}B^{i+j,-j}
\lra 
\bigoplus_{j\in\Bbb Z}B^{i+j+1,-j}
$$
is $d^i_{tot}=\sum_j ((-1)^j d_1^{i+j,-j}+d_2^{i+j,-j})$.
\endgraf
Then the data of a morphism of complexes
$\Phi$ as in i) is by definition the data of a family of maps
$$
\Phi^i: A^i\lra (B^{\point\point})_{tot}^i
$$
such taht the diagram
$$
\xymatrix{
A^i\ar[r]^{d^i_A} \ar[d]_{\Phi^i}& A^{i+1}\ar[d]^{\Phi^{i+1}}\\
(B^{\point\point})_{tot}^i\ar[r]^{d^i_{tot}}&
(B^{\point\point})_{tot}^{i+1}\\
}\leqno{\bf (\NNN)}\rif{comdiff}$$
commutes.
\endgraf
Now $\Phi^i: A^i\lra\bigoplus_j B^{i+j,-j}$, and $B^{\point\point}$ is
bounded, so
$\Phi^i=\sum_j\Phi^i_j$ (which is a finite sum) with
$\Phi^i_j: A^i\lra B^{i+j,-j}$.
The commutativity of (\cite{comdiff}) (regarded in the component 
$B^{i+j+1,-j}$ of $(B^{\point\point})_{tot}^{i+1}$) implies that:
$$
\Phi^{i+1}_j\opp d^i_A= (-1)^j d^{i+j,-j}_1\opp \Phi^i_j
+d^{i+j+1,-(j+1)}_2\opp \Phi^i_{j+1}
$$
so we have the thesis.
}
\smallskip

We use the usual notation $ CF^+_1(\c O_X,\Diff_X)$, resp. $CF^-_1(\c
O_X,\Diff_X)$, 
resp. $CF^b_1(\c O_X,\Diff_X)$ for the category of bounded below, resp. bounded
above, resp. bounded complexes.
We denote by $DF^\ast_1(\c O_X,\Diff_X)$
(resp. $DF^\ast_{1,qc}(\c O_X,\Diff_X)$, resp. 
$DF^\ast_{1,c}(\c O_X,\Diff_X)$ ) 
the full subcategory of $DF_1(\c O_X,\Diff_X)$
whose objects $\c K^{\point}$ satisfy the condition $\Gr^n_F(\c
K^{\point})\in D^{\ast}(\c O_X)$ (resp. $\in D^{\ast}_{\qcoh}(\c O_X)$,
resp. $\in D^{\ast}_{\coh}(\c O_X))$ 
for each $n\in \Bbb N$ and $\ast\in\{+,-,b\}$.